\documentclass[11pt, oneside]{article}   	
\usepackage{geometry}                		
\usepackage{graphicx}				
\usepackage{amsmath,amssymb,amsfonts,graphicx}
\usepackage{multirow,color}
\usepackage{amsmath,amsfonts,amssymb,amsthm, bm}

\usepackage{lipsum}
\usepackage{amsfonts}
\usepackage{graphicx}
\usepackage{epstopdf}
\usepackage{xcolor}
\usepackage{epsfig}

\def\e{\varepsilon}

\newtheorem{remark}{Remark}[section]
\newtheorem{lemma}[remark]{Lemma}
\newtheorem{proposition}[remark]{Proposition}

\newtheorem{theorem}[remark]{Theorem}

\begin{document}

\title{Homogenization in perforated domains at the critical scale}
\author{Giuseppe Cosma Brusca \\ SISSA, Via Bonomea 265, Trieste, Italy \\
gbrusca@sissa.it}
\date{}
\maketitle

\begin{abstract}
    We describe the asymptotic behaviour of the minimal heterogeneous $d$-capacity of a small set, which we assume to be a ball for simplicity, in a fixed bounded open set $\Omega\subseteq \mathbb{R}^d$, with $d\geq2$. Two parameters are involved: $\e$, the radius of the ball, and $\delta$, the length scale of the heterogeneity of the medium. We prove that this capacity behaves as $C|\log \e|^{d-1}$, where $C=C(\lambda)$ is an explicit constant depending on the parameter $\lambda:=\lim_{\e\to0}|\log \delta|/|\log\e|$.
  
    Applying this result, we determine the $\Gamma$-limit of oscillating integral functionals subjected to Dirichlet boundary conditions on periodically perforated domains. In this instance, our first result is used to study the behaviour of the functionals near the perforations which are exactly balls of radius $\e$. 
    We prove that, as in the  homogeneous case, these lead to an additional term that involves $C(\lambda)$.
\end{abstract}


\textbf{\qquad AMS Class:} 49J45, 35B27, 31A15.

\let\thefootnote\relax\footnotetext{Preprint SISSA 02/2023/MATE}

\section{Introduction}

A prototypical variational problem in Sobolev spaces involving scaling-invariant functionals concerns the $d$-capacity of a set $E$ contained in a fixed bounded open set $\Omega\subseteq \mathbb{R}^d$ with $d\geq2$. If we assume $E$ to have diameter of size $\e\ll1$, an explicit computation proves that the asymptotic behaviour of such capacity equals $|\log\e|^{1-d}$, up to a dimensional factor.

In this paper, we introduce a dependence on $x$, which in the model describes the heterogeneity of a medium, and we analyse the asymptotic behaviour of minima 
\begin{equation}\label{minpbm1}
m_{\e,\delta}:=\min\Bigl\{\int_{\Omega} f\left(\frac{x}{\delta},\nabla u(x)\right)\,dx: u\in W^{1,d}_0(\Omega), u=1 \hbox{ on } B(z,\e), z\in\Omega\Bigr\},
\end{equation} 
where $\delta=\delta(\e)$ is positive and vanishing as $\e\to0$, and $f(x,\xi)$ is a function with suitable assumptions of periodicity and homogeneity.

\noindent We assume $f:\mathbb{R}^d \times \mathbb{R}^d \to [0,+\infty)$ to be a Borel function with the following properties:

\smallbreak

(P$1$) (periodicity) $f(\cdot,\xi) \hbox{ is } 1\hbox{-periodic for every } \xi\in\mathbb{R}^d$, i.e., denoting by $e_k$ an element of the canonical basis, $f(x+e_k,\xi)=f(x,\xi)$ for every $x$ and $\xi$ in $\mathbb{R}^d$, and $k=1,...,d\,$;

\smallbreak

(P$2$) (positive $d$-homogeneity) $f(x,t\,\xi)=t^d f(x,\xi) \hbox{ for every } t>0$, for every  $x$ and $\xi$ in $\mathbb{R}^d$;

\smallbreak

(P$3$) (standard growth conditions of order $d$) there exist $\alpha, \beta$ such that $0<\alpha <\beta$ and $\alpha|\xi|^d\leq f(x,\xi)\leq \beta|\xi|^d$ for every $x$ and $\xi$ in $\mathbb{R}^d$.

\medbreak

In light of (P$1$) and (P$2$), the minima defined in \eqref{minpbm1} stand for the minimal \emph{heterogeneous capacity} of a small set (which is not restrictive to assume to be a ball) of size $\e$, while $\delta$ is the period of the heterogeneity modelled by oscillating terms.

\noindent Assumption (P$3$) is technical as it is needed to apply the Homogenization Theorem. 

We remark that, by a relaxation argument, we may assumw $f(x,\xi)$ being convex on the second variable so that the associated energy functional is $W^{1,d}(\Omega)$-weakly lower semicontinuous and the terms defined by \eqref{minpbm1} are actual minima. 

\medbreak

The first result we achieve is the asymptotic estimate for the minima in \eqref{minpbm1}. To this end, we work along subsequences (not relabeled) for which it exists
\[
\lambda := \lim_{\e\to0}\frac{|\log\delta|}{|\log\e|} \land 1 \in [0,1].
\]
We introduce a function describing the asymptotic concentration of the heterogeneous capacity at a point $z\in\mathbb{R}^d$; it is given by
\begin{equation}\label{Phi}
\begin{split}
\Phi(z):=\lim_{R\to+\infty}(\log R)^{d-1}\min\Bigl\{\int_{B(0,R)\setminus B(0,1)} f(z,\nabla u(x))\,dx: \, & u\in W^{1,d}_0(B(0,R)), \\
& u=1 \hbox{ on } B(0,1)\Bigr\},
\end{split}
\end{equation}
then we define a constant portraying the effect of homogenization, which is
\begin{equation}\label{Chom}
\begin{split}
C_{\rm hom}:= \lim_{R\to+\infty} (\log R)^{d-1}\min\Bigl\{\int_{B(0,R)\setminus B(0,1)} f_{\rm hom}(\nabla u(x))\,dx : \, & u\in W^{1,d}_0(B(0,R)), \\
& u=1 \hbox{ on } B(0,1)\Bigr\},
\end{split}
\end{equation}
where $f_{\rm hom}$ is the positively $d$-homogeneous function 
\begin{equation}\label{fhom}
f_{\rm hom}(\xi)= \min \Bigl\{ \int_{(0,1)^d} f(y, \xi+ \nabla \varphi(y))\,dy : \varphi \in W^{1,d}_{loc}(\mathbb{R}^d), \varphi \hbox{ 1-periodic} \Bigr\}
\end{equation}
determined by the Homogenization Theorem.

We prove that if there exists a point $x_0\in\mathbb{R}^d$ for which $\Phi$ is consistently concentrated in an optimal way, then it holds
\begin{equation}\label{result1}
\lim_{\e\to0}|\log\e|^{d-1}m_{\e,\delta} = \Phi(x_0)C_{\rm hom}\Bigl[\lambda\Phi(x_0)^{\frac{1}{d-1}}+(1-\lambda)C_{\rm hom}^{\frac{1}{d-1}}\Bigr]^{1-d}. 
\end{equation}

\noindent As an example, we refer to the quadratic case, already treated in \cite{BraBru}. 
If $d=2$ and $f(x,\xi)=a(x)|\xi|^2$, where $a(x)$ is a $1$-periodic continuous function bounded form below by a constant $\alpha$, then we can pick $x_0$ so that $\Phi(x_0)=2\pi\alpha$, and denoting the homogenized matrix by $A_{\rm hom}$, we obtain $C_{\rm hom} = 2\pi\sqrt{\det A_{\text{\rm hom}}}$. 
We eventually find
\[
\lim_{\e\to0}|\log\e|m_{\e, \delta} = 2\pi\frac{\alpha\sqrt{\det A_{\text{\r hom}  }}}{\lambda\alpha+(1-\lambda)\sqrt{\det A_{\text{\r hom}  }}}\,.
\]

\medbreak

A fundamental tool in the proof of this result is a method elaborated by De Giorgi, which allows to impose boundary conditions on converging sequences. In this work, it is presented and proved in a version (Lemma \ref{modlemma}) which is suitable to our purposes, similar to that in \cite{AnsBra}. 

\bigbreak

The second result concerns homogenization on perforated domains. The argument follows the work of Laura Sigalotti in \cite{Sigalotti}; that is, a nonlinear version at the critical exponent of the homogenization on perforated domains with Dirichlet boundary conditions originally studied, e.g., by Marchenko and Khruslov in \cite{MK} and by Cioranescu and Murat in \cite{CM}. Works about the asymptotic behaviour of Dirichlet problems in varying domains are, e.g., \cite{BL, CDG, DMM}, or also \cite{LLL} for the numerical perspective.  

\smallbreak

Denoting by $B$ the open unit ball, and by $d(\e)$ the period of the perforations, we define a periodically perforated domain as
\[
\Omega_{\e}:=\Omega \setminus \bigcup_{i\in \mathbb{Z}^d} id(\e)+\e B
\]
and we consider functionals $F_{\e}: L^d(\Omega)\rightarrow[0,+\infty]$ given by
\[
F_{\e}(u):=
\begin{cases}
\displaystyle\int_{\Omega}f\left(\frac{x}{\delta}, \nabla u(x)\right)\,dx & \hbox{ if } u\in W^{1,d}(\Omega) \hbox{ and } u=0 \hbox{ on } \Omega\setminus\Omega_{\e} \\[5pt]
+\infty & \hbox{ otherwise}.
\end{cases}
\]

In the above mentioned works \cite{CM} and \cite{MK}, it is analysed the homogeneous case $f(x,\xi)=|\xi|^p$ for $p>1$, and it is provided a critical choice for the period, which is exactly $d(\e)=|\log \e|^{\frac{1-d}{d}}$, if $p=d$. Moreover, the $\Gamma$-limit with respect to the strong convergence in $L^d(\Omega)$ is proved to be 
\[
\int_{\Omega}|\nabla u(x)|^d\,dx + \kappa_d \int_{\Omega}|u(x)|^d\,dx
\]
for every $u\in W^{1,d}(\Omega)$, with $\kappa_d$ a dimensional constant, showing that internal boundary conditions disappear with the arising of a so-called \emph{strange term}.

\smallbreak

We prove an analogous statement: for simplicity, we assume that $d(\e)$ is an integer multiple of $\delta(\e)$ so that the periodicity of the perforation is 'compatible' with that of the energy. As an oscillating term is introduced in our study, we expect our result to be affected by the different rate of vanishing of $\e,\delta$ and since this information is encoded by the parameter $\lambda$, we show that for every $u\in W^{1,d}(\Omega)$ it holds 
\[
\Gamma\emph{-}\lim_{\e} F_\e(u)=\int_{\Omega}f_{\rm hom}(\nabla u(x))\,dx + C(\lambda)\int_{\Omega}|u(x)|^d\,dx,
\]
where $C(\lambda)$ is the constant
\[
C(\lambda):=
\Phi(0)C_{\rm hom}\Bigl[\lambda\Phi(0)^{\frac{1}{d-1}}+(1-\lambda)C_{\rm hom}^{\frac{1}{d-1}}\Bigr]^{1-d},
\]
and the term $\Phi(0)$ is due to the asymptotic analysis of the problems
\begin{equation*}
\min\Bigl\{\int_{B} f\left(\frac{x}{\delta},\nabla u(x)\right)\,dx: u\in W^{1,d}_0(B), u=1 \hbox{ on } B(0,\e)\Bigr\},
\end{equation*} 
where the centres of the perforations have been fixed at $0$ for every $\e$.

\subsection{Preliminaries}

In this section and the following ones, let $d\geq2$, $\Omega\subseteq \mathbb{R}^d$ be a bounded open set and $\lambda :=\lim_{\e\to0}|\log\delta|/|\log\e|$.

We start by justifying the definitions given in \eqref{Phi} and \eqref{Chom} thorugh the following lemma which takes advantage of a scaling invariance argument.

\begin{lemma}\label{limlemma}
Let $g: \mathbb{R}^d\rightarrow \mathbb{R}$ be a Borel function which is positively homogeneous of degree $d$ and assume there exist positive constants $C_1< C_2$ such that $C_1|\xi|^d\leq g(\xi)\leq C_2|\xi|^d$ for every $\xi\in\mathbb{R}^d$. Define 
\[
m_{R}:=\min\Bigl\{\int_{B(0,R)\setminus B(0,1)} g(\nabla u(x))\,dx: u\in W^{1,d}_0(B(0,R)), u=1 \hbox{ on } B(0,1)\Bigr\},
\]
then it exists $\lim\limits_{R\to+\infty} (\log R)^{d-1}\,m_{R}$ and this limit is finite.
\end{lemma} 

\begin{proof}
Fix $S>R$ and put $T:=\lfloor \log S/\log R\rfloor$ so that the annuli $B(0,R^k)\setminus B(0,R^{k-1})$ are contained in $B(0,S)\setminus B(0,1)$ for every $k=1,...,T$.

Let $u$ be a solution of the problem 
\[
\min\Bigl\{\int_{B(0,R)\setminus B(0,1)} g(\nabla u(x))\,dx : u\in W^{1,d}_0(B(0,R)), u=1 \hbox{ on } B(0,1)\Bigr\},    
\]
and for $k=1,...,T$, define functions $u^k\in W^{1,d}(B(0,R^k)\setminus \overline{B}(0,R^{k-1}))$ as 
\[
u^k(x):=\frac{1}{T}u\left(\frac{x}{R^{k-1}}\right)+\frac{T-k}{T};
\]
then put
$u_S\in W^{1,d}_0(B(0,S))$ as
\[
u_S(x):= 
\begin{cases}
1 & \hbox{ if } x\in B(0,1) \\
u^k(x) & \hbox{ if } x\in  B(0,R^k)\setminus B(0,R^{k-1}),\, k=1,...,T \\
0 & \hbox{ if } x\in B(0,S)\setminus B(0,R^T).
\end{cases}
\]

We have
\begin{eqnarray} \nonumber 
(\log S)^{d-1}m_{S} &\leq& (\log S)^{d-1} \int_{B(0,S)\setminus B(0,1)}  g(\nabla u_S(x))\,dx \\ \nonumber
&=& (\log S)^{d-1} \sum_{k=1}^T \int_{B(0,R^k)\setminus B(0,R^{k-1})} g(\nabla u^k(x))\,dx \\ \nonumber
&=& (\log S)^{d-1} \sum_{k=1}^T \frac{1}{T^d}\int_{B(0,R)\setminus B(0,1)} g(\nabla u(x))\,dx \\ \nonumber
&=& (\log S)^{d-1} \frac{1}{T^{d-1}}m_R \\ \nonumber
&\leq& (\log S)^{d-1} \Bigl(\frac{\log R}{\log S-\log R}\Bigr)^{d-1}m_{R}.
\end{eqnarray}

If we pass to the $\limsup$ as $S\to+\infty$, and then we pass to the $\liminf$ as $R\to+\infty$, we obtain
\[
\limsup_{S\to+\infty}\, (\log S)^{d-1}m_{S} \leq \liminf_{R\to+\infty}\, (\log R)^{d-1}m_{R}.
\]
 
In order to check the limit is finite, consider the function
\[
u(x):=1-\frac{\log |x|}{\log R},\qquad x\in B(0,R)\setminus \overline{B}(0,1),
\]
and note that the estimate
\begin{multline*}
(\log R)^{d-1}m_{R} \leq \int_{B(0,R)\setminus B(0,1)} g(\nabla u(x)) = (\log R)^{d-1} \int_{B(0,R)\setminus B(0,1)} g\biggl(\frac{x}{-|x|^2\log R}\biggr)\,dx \\
= (\log R)^{-1} \int_{B(0,R)\setminus B(0,1)} g\biggl(-\frac{x}{|x|^2}\biggr)\,dx  \leq (\log R)^{-1} C_2\int_{B(0,R)\setminus B(0,1)} \frac{1}{|x|^d}\,dx = C_2 \sigma_{d-1}
\end{multline*}
holds, completing the proof.
\end{proof}

\medbreak

We state the Homogenization Theorem (see \cite{GCB, BDF, DM, JKO}) in a slightly modified version which allows to take into account also translations in the following sense. 

\begin{theorem}\label{homthmtranslation}
Let $A$ be a bounded open subset of $\mathbb{R}^d$ with Lipschitz boundary and $(\tau_{\eta})_{\eta>0}\subseteq \mathbb{R}^d$. Then 
\[
\Gamma\emph{-}\lim_{\eta\to0} \int_{A}f\left(\frac{x}{\eta}+\tau_{\eta}, \nabla u(x)\right)\,dx = \int_{A}f_{\rm hom}(\nabla u(x))\,dx\,,
\]
for every $u\in W^{1,d}(A)$, where the $\Gamma$-limit is meant with respect to the strong convergence in $L^d(A)$ and $f_{\rm hom}$ is the function given by \eqref{fhom}.

In particular, for every $\phi\in W^{1,d}(A)$ we have
\begin{multline*}
\lim_{\eta\to0}\inf\Bigl\{\int_{A}f\left(\frac{x}{\eta}+\tau_{\eta}, \nabla u(x)\right)\,dx : u\in\phi+W^{1,d}_0(A)\Bigr\} \\
= \min\Bigl\{\int_{A}f_{\rm hom}(\nabla u(x))\,dx : u\in\phi+W^{1,d}_0(A)\Bigr\}\,.
\end{multline*}
\end{theorem}

\medbreak

At this point, Lemma \ref{limlemma} and assumptions (P$2$), (P$3$), make well defined the function
\begin{equation*}
\begin{split}
\Phi(z):=\lim_{R\to+\infty}(\log R)^{d-1}\min\Bigl\{\int_{B(0,R)\setminus B(0,1)} f(z,\nabla u(y))\,dy: \, & u\in W^{1,d}_0(B(0,R)), \\
& u=1 \hbox{ on } B(0,1)\Bigr\},
\end{split}
\end{equation*}
while, to introduce properly the constant
\begin{equation*}
\begin{split}
C_{\rm hom}:= \lim_{R\to+\infty} (\log R)^{d-1}\min\Bigl\{\int_{B(0,R)\setminus B(0,1)} f_{\rm hom}(\nabla u(x))\,dx : \, & u\in W^{1,d}_0(B(0,R)), \\
& u=1 \hbox{ on } B(0,1)\Bigr\},
\end{split}
\end{equation*}
we also need to rely on Theorem \ref{homthmtranslation}: this, combined with the fact that the growth conditions on $f(x,\xi)$ posed in (P$3$) are inherited by the function $f_{\rm hom}(\xi)$, ensures that the above lemma applies.

\section{Asymptotic analysis of minima}

We first aim at estimating the asymptotic behaviour of the minima with fixed centres modulo a translation. More precisely, let $z$ be a point in $\Omega$, for every $\e$ positive sufficiently small, let $(z_\e)_\e$ be a family of points in $\Omega$ of the form $z_\e=\delta z+\delta i_\e$, where $(i_\e)_\e\subseteq \mathbb{Z}^d$. Also assume that such family of points is well contained in $\Omega$, i.e., that $\inf_{\e} \text{dist}(z_\e, \partial \Omega)>0$; we put
\begin{equation}\label{minpbm2}
\mu_{\e,\delta}=\min\Bigl\{\int_{\Omega} f\left(\frac{x}{\delta},\nabla u(x)\right)\,dx: u\in W^{1,d}_0(\Omega), u=1 \hbox{ on } B(z_{\e},\e)\Bigr\}.
\end{equation}

The asymptotic behaviour of these minima is the main subject of this section; we prove the following.

\begin{proposition}\label{mainpropfixed} 
Let $z\in \Omega$ be a fixed point, and let $(z_{\e})_{\e}$ be a family of points equal to $z$ modulo $\delta$ as above. Assume that for every $\nu>0$, there exists $r_{\nu}>0$ such that for every $x\in B(z, r_{\nu})$ it holds 
\begin{equation}\label{assumption0}
|f(z,\xi)-f(x,\xi)|\leq \nu|\xi|^d\, \hbox{ for every } \xi\in\mathbb{R}^d.
\end{equation}
\noindent Then 
\[
\lim_{\e\to0}|\log\e|^{d-1}\mu_{\e,\delta}=
\Phi(z)C_{\rm hom}\Bigl[{\lambda\Phi(z)^{\frac{1}{d-1}}+(1-\lambda)C_{\rm hom}^{\frac{1}{d-1}}}\Bigr]^{1-d}.
\]
\end{proposition}

\bigbreak

The proof is divided in two parts: the bound from below and the construction of an optimal sequence. In the first, the main tool we use is the following lemma which allows to modify a function in order to attain constant values (in the sense of the trace) on the boundary of a thin annulus, still controlling the value of the associated energy.

\begin{lemma}\label{modlemma}
Let $f: \mathbb{R}^d \times \mathbb{R}^d\rightarrow \mathbb{R}$ be a Borel function satisfying the standard growth conditions property {\rm (P$3$)}. Let $z\in\mathbb{R}^d$, $R>0$ and define
\[
F(u,A):=\,\int_{A} f(x,\nabla u(x))\,dx
\]
for every $u\in W^{1,d}(B(z,R))$ and $A\subseteq B(z,R)$ Borel subset.

Let $\eta>0$, put $S:=\max\,\{s\in \mathbb{N} : \eta2^s \leq R\}$ and assume $S\geq3$. Take $N$ natural number such that $2\leq N<S$ and $r$ positive real number such that $r \leq \eta2^{S-N}$.

\noindent Then there exists a function $v$ with the following properties:

\emph{(i)} $v\in W^{1,d}(B(z,R)\setminus \overline{B}(z,r))$, 

\emph{(ii)} $\hbox{there exists } \,j\in\{1,...,N-1\} \hbox{ such that }$
\[
v=u \hbox{ on } B(z,\eta2^{S-j-1})\setminus \overline{B}(z,r))\cup B(z,R)\setminus \overline{B}(z,\eta2^{S-j+1}),
\] 

\emph{(iii)}$\hbox{ for the same } j, \hbox{ the function }  v \hbox{ is constant on } \partial B(z,\eta2^{S-j})$, 

\emph{(iv)} There exists a positive constant $C$ depending on $\alpha, \beta$ and the dimension $d$ such that
\[
F(v,B(z,R)\setminus B(z,r))\leq \Bigr(1+\frac{C}{N-1}\Bigr)F(u,B(z,R)\setminus B(z,r)).
\]
\end{lemma}
\medbreak

\begin{proof} Assume $z=0$, if not, center the construction around $z$ and repeat the argument.

For $k=1,...,N-1$, we define annuli $A_k:=B(0,\eta 2^{S-N+k+1})\setminus B(0,\eta 2^{S-N+k-1})$ and radial cutoff functions
\[
\phi_k(\rho):=
\begin{cases}
0 & \hbox{ if } \rho\in[0,\eta2^{S-N+k-1}] \\[5pt]
\frac{\rho-\eta2^{S-N+k-1}}{\eta2^{S-N+k-1}} & \hbox{ if } \rho\in(\eta2^{S-N+k-1}, \eta2^{S-N+k}] \\[5pt]
\frac{\eta2^{S-N+k+1}-\rho}{\eta2^{S-N+k}} & \hbox{ if } \rho\in(\eta2^{S-N+k},\eta2^{S-N+k+1}] \\[5pt]
0 & \hbox{ if } \rho\in(\eta2^{S-N+k+1},R],
\end{cases}
\]
then we put $\psi_k:=1-\phi_k$ and define $v_k:=\psi_ku+(1-\psi_k)u_{A_k}$, where we denote by $u_{A_k}$ the integral average of $u$ on $A_k$. 

At each fixed $k$, taking into account that $|\psi_k| \leq 1$ and 
\[
|\nabla \psi_k|^d=|\nabla \phi_k|^d \leq \Bigl(\frac{1}{\eta2^{S-N+k-1}}\Bigr)^d,
\]
we exploit (P$2$) to have
\begin{equation}\label{pw}
\begin{split}
\int_{A_k} f(x,\nabla v_k(x))\,dx & \leq  \beta\int_{A_k} |\nabla v_k(x)|^d\,dx \\
& = \beta\int_{A_k} |\psi_k\nabla u(x)+(u-u_{A_k})\nabla\psi_k(x)|^d\,dx
\\
& \leq \beta2^{d-1} \Bigl[\int_{A_k} |\nabla u|^d\,dx  +\Bigl(\frac{1}{\eta 2^{S-N+k-1}}\Bigr)^d\int_{A_k} |u(x)-u_{A_k}|^d\,dx \Bigr].
\end{split}
\end{equation}
Consider now the following well known scaling property of the Poincaré-Wirtinger inequality: given $A$ open, bounded, connected, with Lipschitz boundary and $\lambda>0$, it holds
\[
\frac{1}{\lambda ^d}\int_{\lambda A} |u-u_{\lambda A}|^d\,dx \leq P(A) \int_{\lambda A}|\nabla u|^d\,dx,
\]
where $u_{\lambda A}$ is the integral average of $u$ on $\lambda A$ and $P(A)$ is the Poincaré-Wirtinger constant related to $A$.

\noindent We apply this result with $A=B(0,4) \setminus \overline{B}(0,1)$ and $\lambda=\eta 2^{S-N+k-1}$, obtaining
\[
\Bigl(\frac{1}{\eta 2^{S-N+k-1}}\Bigr)^d\int_{A_k} |u(x)-u_{A_k}|^d\,dx  \leq
P^d \int_{A_k} |\nabla u|^d\,dx,
\]
being $P:=P(A)$ a constant which does not depend on $k$.

\noindent As a consequence \eqref{pw} turns into
\[
\begin{split}
\int_{A_k} f(x,\nabla v_k(x))\,dx & \leq \beta2^{d-1}\bigl(1+P^d\bigr)\int_{A_k} |\nabla u|^d\,dx \\
& \leq \frac{\beta}{\alpha}2^{d-1}\bigl(1+P^d\bigr)\int_{A_k} f(x,\nabla u(x))\,dx,
\end{split}
\]
and summing over $k$, we deduce
\[
\sum_{k=1}^{N-1} \int_{A_k} f(x,\nabla v_k(x))\,dx \leq C\int_{B(0,R)\setminus B(0,r)} f(x,\nabla u(x))\,dx,
\]
where we put $C:=\beta2^{d-1}\bigl(1+P^d\bigr)/\alpha$. It follows there exists $j\in\{1,...,N-1\}$ such that
\[
\int_{A_j} f(x,\nabla v_j(x))\,dx \leq \frac{C}{N-1}\int_{B(0,R)\setminus B(0,r)} f(x,\nabla u(x))\,dx,
\]
and then it holds
\[
\begin{split}
\int_{B(0,R)\setminus B(0,r)} f(x,\nabla v_j(x))\,dx & = \int_{(B(0,R)\setminus B(0,r))\setminus A_j} f(x,\nabla u(x))\,dx + \int_{A_j} f(x,\nabla v_j(x))\,dx \\
& \leq \left(1+\frac{C}{N-1}\right)\int_{B(0,R)\setminus B(0,r)} f(x,\nabla u(x))\,dx,
\end{split}
\]
which concludes the proof.
\end{proof}

The estimate in (iv) is more precise as $N\to\infty$, i.e., as $\eta\to0$. Our strategy will consist in parting the open set $\Omega$ through many annuli having small inner and outer radii, say of order $\e^{\lambda}\sim \delta$, and there modifying a function $u\in W^{1,d}_0(\Omega)$ to achieve some constant Dirichlet boundary conditions as a consequence of (iii). The error introduced by the modification will be negligible in light of (iv).  

\subsection{Lower bound}

In what follows, we systematically identify a function $u\in W^{1,d}_0(\Omega)$ with the the extension obtained by setting $u=0$ on $\mathbb{R}^d\setminus\Omega$, which belongs to $W^{1,d}(\mathbb{R}^d)$.

\noindent For simplicity of notation, given $A$ a Borel subset of $\mathbb{R}^d$ and $u\in W^{1,d}(\mathbb{R}^d)$, we put
\[
F_{\e}(u, A):= \int_A f\left(\frac{x}{\delta}, \nabla u(x)\right)\,dx
\]
and denote by $R_{\Omega}$ the maximum among the diameter of $\Omega$ and $1$.

\smallbreak

We consider separately the cases $\lambda=0,\, \lambda\in(0,1)$ and $\lambda=1$; we obtain for each instance the same kind of estimate and then we conclude by the same argument.

\smallbreak

If $\lambda=0$, fix a parameter $\lambda_2\in(\lambda,1)$ so that
\[
\frac{\e^{\lambda_2}}{\delta}\to0 \hbox{ as } \e\to0.
\]
For every $u\in W^{1,d}_0(\Omega)$ such that $u=1$ on $B(z_{\e},\e)$, the inclusion $\Omega \subseteq B(z_{\e}, R_{\Omega})$ leads to the equality
\[
F_{\e}(u,\Omega) = F_{\e}(u,B(z_{\e}, R_{\Omega})),
\]
then we apply Lemma \ref{modlemma} to the function $u\in W^{1,d}_0(B(z_\e, R_\Omega))$, with 
\[
f(x,\xi)=f\left(\frac{x}{\delta},\xi\right),\, \eta=\e,\, R=\e^{\lambda_2}\,, N \in \mathbb{N}\cap \left(1, \left\lfloor \frac{(1-\lambda_2)|\log\e|}{\log2} \right\rfloor=S \right) \,\hbox{ and }  r=\e.
\]
We get a function $v\in W^{1,d}_0(B(z_{\e}, R_{\Omega}))$ such that $v=1$ on $B(z_{\e},\e)$, $v=c$ on $\partial B(z_{\e}, \e2^{S-j})$ for some constant $c$ and some index $j\in\{1,...,N-1\}$, and $v=u$ on $B(z_{\e}, R_{\Omega})\setminus B(z_{\e},\e^{\lambda_2})$; hence, it holds
\begin{equation}\label{below1lambda0}
\begin{split}
\left(1+\frac{C}{N-1}\right) F_\e(u,\Omega) & = \left(1+\frac{C}{N-1}\right)F_{\e}(u,B(z_{\e}, R_{\Omega}))  \geq F_{\e}(v,B(z_{\e}, R_{\Omega})) \\
& = F_{\e}(v,B(z_{\e},\e2^{S-j})) + F_{\e}(v,B(z_{\e}, R_{\Omega})\setminus B(z_{\e},\e2^{S-j})).
\end{split}
\end{equation}
Now we set 
\[
w^1:=
\begin{cases}
v & \hbox{ on } B(z_{\e}, \e2^{S-j}) \\
c & \hbox{ on } B(z_{\e}, \e^{\lambda_2})\setminus B(z_{\e}, \e2^{S-j})
\end{cases}\qquad
w^2:=
\begin{cases}
c & \hbox{ on } B(z_{\e}, \e2^{S-j})\setminus \overline{B}(z_{\e}, \e2^{S-N}) \\
v & \hbox{ on } B(z_{\e}, R_{\Omega})\setminus  B(z_{\e}, \e2^{S-j}),
\end{cases}
\]
and we note that
\[
F_{\e}(w^1,B(z_{\e},\e^{\lambda_2})) = F_{\e}(v,B(z_{\e},\e2^{S-j}))
\]
and
\[
F_{\e}(w^2,B(z_{\e}, R_{\Omega})\setminus B(z_{\e},\e2^{S-N})) =
F_{\e}(v,B(z_{\e}, R_{\Omega})\setminus B(z_{\e},\e2^{S-j})),  
\]
thus 
\[
F_{\e}(v,B(z_{\e}, R_{\Omega})) = F_{\e}(w^1,B(z_{\e}, \e^{\lambda_2})) +  F_{\e}(w^2, B(z_{\e}, R_{\Omega})\setminus  B(z_{\e}, \e2^{S-N})).
\]

At this point we take advantage of the fact that both the functions $w^1$ and $w^2$ attain constant values on the components of the boundary of their domain. 

\noindent We rewrite inequality \eqref{below1lambda0} as
\[
\left(1+\frac{C}{N-1}\right)F_\e(u,\Omega) \geq 
\]
\begin{equation*}
\begin{split}
& \geq \min\{F_{\e}(v,B(z_{\e},\e^{\lambda_2})) : v\in W^{1,d}(B(z_{\e},\e^{\lambda_2})), v=1 \hbox{ on } B(z_{\e},\e), v=c \hbox{ on } \partial B(z_{\e},\e^{\lambda_2})\} \\
& + \min\{F_{\e}(v,B(z_{\e}, R_{\Omega})\setminus \overline{B}(z_{\e},\e2^{S-N})) : 
v\in W^{1,d}(B(z_{\e}, R_{\Omega})\setminus \overline{B}(z_{\e},\e^{\lambda_2}2^{S-N})),    \\
& \qquad \qquad \qquad \qquad \qquad \qquad v=c \hbox{ on } \partial B(z_{\e},\e2^{S-N}), v=0 \hbox{ on } \partial B(z_{\e}, R_{\Omega})\} ,
\end{split}
\end{equation*}
and taking into account the transformations 
\[
v(x) \mapsto \frac{v(x)-c}{1-c}\,, \qquad  \qquad v(x)\mapsto \frac{v(x)}{c}\,,
\] 
and the property of homogeneity (P$2$), we have that the last expression equals
\begin{align} \label{min1}
& \min\{F_{\e}(v,B(z_{\e},\e^{\lambda_2})) : v\in W^{1,d}_0(B(z_{\e},\e^{\lambda_2})), 
v=1 \hbox{ on } B(z_{\e},\e)\}|1-c|^d \\ \nonumber
+& \min\{F_{\e}(v,B(z_{\e}, R_{\Omega})\setminus  \overline{B}(z_{\e},\e2^{S-N})):v\in W^{1,d}(B(z_{\e}, R_{\Omega})\setminus \overline{B}(z_{\e},\e^{\lambda_2}2^{S-N})), \\ \label{min2}
& \qquad \qquad \qquad \qquad \qquad \qquad v=1 \hbox{ on } B(z_{\e},\e2^{S-N}), v=0 \hbox{ on } \partial B(z_{\e}, R_{\Omega})\}|c|^d.
\end{align}

\smallskip

To treat the minimum in \eqref{min1}, we apply the transformation $v(x)\mapsto v(z_{\e}+x\e)$, getting
\[
\min\Bigl\{\int_{B(0,\e^{\lambda_2-1})} f\Bigl(\frac{x}{\delta}\e+\frac{z_{\e}}{\delta},\nabla v(x)\Bigr) : v\in W^{1,d}_0(B(0,\e^{\lambda_2-1})), v=1 \hbox{ on } B(0,1) \Bigr\}|1-c|^d\,.
\]
As $z_{\e}=\delta z+\delta i_{\e}$, we exploit the periodicity assumption (P$1$) to get 
\[
f\left(\frac{x}{\delta}\e+\frac{z_{\e}}{\delta},\xi\right)=f\left(\frac{x}{\delta}\e+z,\xi\right) \hbox{ for every } \xi\in\mathbb{R}^d;
\]
also note that if $x\in B(0,\e^{{\lambda_2}-1})$, then $\frac{\e}{\delta}|x|<\frac{\e^{\lambda_2}}{\delta} \to0$ as $\e\to0$. Hence, for every $\nu>0$, given $r_\nu$ as in \eqref{assumption0}, it holds that $B(0,\e^{\lambda_2-1}) \subseteq B(0, r_\nu)$, so that for every $\e$ sufficiently small we have
\[
f\Bigl(\frac{x}{\delta}\e+\frac{z_{\e}}{\delta},\xi\Bigr) \geq f(z, \xi)-\nu|\xi|^d \hbox{ for every } \xi\in\mathbb{R}^d.
\]
Combining these observations with the growth condition (from below) in (P$3$), we get
\[
\int_{B(0,\e^{\lambda_2-1})} f\Bigl(\frac{x}{\delta}\e+\frac{z_{\e}}{\delta},\nabla v(x)\Bigr)\,dx \geq \left(1-\frac{\nu}{\alpha}\right) \int_{B(0,\e^{\lambda_2-1})} f(z,\nabla v(x))\,dx
\]
for every $v\in W^{1,d}_0(B(0,\e^{\lambda_2-1}))$ such that $v=1$ on $B(0,1)$. By the application of Lemma \ref{limlemma}, which is possible due to the fact that $\lambda_2<1$, we obtain
\[ 
\min\{F_{\e}(v,B(z_{\e},\e^{\lambda_2})) : v\in W^{1,d}_0(B(z_{\e},\e^{\lambda_2})), 
v=1 \hbox{ on } B(z_{\e},\e)\}|1-c|^d 
\]
\begin{equation}\label{below2lambda0}
\geq
\frac{\Phi(z)+o_{\e}(1)}{(1-\lambda_2)^{d-1}|\log\e|^{d-1}}|1-c|^d,
\end{equation}
where we get rid of the term in $\nu$ since $1-\nu/\alpha$ may be taken arbitrarily close to $1$ as $\e\to0$.

\smallbreak

In order to deal with the minimum in \eqref{min2}, we apply once more property (P$3$), and in particular the inequality $f(x, \xi) \geq \alpha |\xi|^d$. We get a lower bound in terms of the d-capacity of the inclusion $B(z_{\e},\e2^{S-N})\subseteq B(z_{\e}, R_{\Omega})$ which is explicitly computed; more precisely, we have
\begin{multline*}
\min\{F_{\e}(v,B(z_{\e}, R_{\Omega})\setminus\overline{B}(z_{\e},\e2^{S-N})):v\in W^{1,d}(B(z_{\e}, R_{\Omega})\setminus \overline{B}(z_{\e},\e^{\lambda_2}2^{S-N})), \\
v=1 \hbox{ on } B(z_{\e},\e2^{S-N}), v=0 \hbox{ on } \partial B(z_{\e}, R_{\Omega})\}|c|^d 
\end{multline*}
\begin{eqnarray} \notag
&\geq& \alpha \hbox{Cap}_d(B(z_{\e},\e2^{S-N}),B(z_{\e}, R_{\Omega}))|c|^d \\ \notag
& = & \frac{\alpha\sigma_{d-1}}{[\log R_{\Omega}+|\log\e|-(S-N)\log2]^{d-1}}|c|^d \\ \label{below3lambda0}
& \geq & \frac{\alpha\sigma_{d-1}}{[\log  R_{\Omega}+\lambda_2|\log\e|+(N+2)\log2]^{d-1}}|c|^d,
\end{eqnarray}
where the last inequality follows recalling that $S=\lfloor \frac{(1-\lambda_2)|\log\e|}{\log2} \rfloor$.

\smallbreak

Gathering \eqref{below2lambda0} and \eqref{below3lambda0}, and multiplying by $|\log\e|^{d-1}$, we get
\begin{equation}\label{below4lambda0}
\begin{split}
\left(1+\frac{C}{N-1}\right)|\log\e|^{d-1}F_{\e}(u,\Omega) & \geq \frac{\Phi(z)+o_{\e}(1)}{(1-\lambda_2)^{d-1}}|1-c|^d \\
& + \frac{\alpha\sigma_{d-1}|\log\e|^{d-1}}{[\log R_{\Omega}+\lambda_2|\log\e|+(N+2)\log2]^{d-1}}|c|^d\,.
\end{split}
\end{equation}

We recall that, by construction, the constant boundary value $c$ actually depends on $\e$, being the mean value of the function $u$ in an annulus whose radii are $\e$-dependent.  In order to pass to the lower limit as $\e\to0$, we make precise the fact that we can assume $c(\e)\to c\in \mathbb{R}$.

\noindent An easy way to see this is observing that we may assume that $u$ takes values in $[0,1]$ as trivially follows by the estimate
\[
F_\e(u, \Omega) \geq F_\e((u \lor 0) \land 1,\Omega) \hbox{ for every } u\in W^{1,d}_0(\Omega),
\]
so that $c(\e)\in[0,1]$ as well.
\noindent Then we find a sequence $\e_k\to0$, and correspondingly $c_k:=c(\e_k)$, such that
\[
\liminf_{k \to+\infty} \frac{\Phi(z)+o_{k}(1)}{(1-\lambda_2)^{d-1}}(1-c_k)^d + \frac{\alpha\sigma_{d-1}|\log\e_k|^{d-1}}{[\log R_{\Omega}+\lambda_2|\log\e_k|+(N+2)\log2]^{d-1}}(c_k)^d
\]
is achieved as a limit; as $(c_k)_k\subseteq [0,1]$, we extract a further subsequence $c_{k_h}\to c\in[0,1]$, to get
\[
\left(1+\frac{C}{N-1}\right) \liminf_{\e\to0} |\log\e|^{d-1}F_{\e}(u,\Omega)  \geq \frac{\Phi(z)}{(1-\lambda_2)^{d-1}}|1-c|^d + \frac{\alpha\sigma_{d-1}}{\lambda_2^{d-1}}|c|^d\,.
\]

\noindent Finally, we pass to the limit as $N\to+\infty$ and recall that $u$ was arbitrary among the admissible functions for the minimization; we conclude that 
\begin{equation}\label{importantliminf1}
\liminf_{\e\to0}|\log\e|^{d-1}\mu_{\e,\delta} \geq \frac{\Phi(z)}{(1-\lambda_2)^{d-1}}|1-c|^d + \frac{\alpha\sigma_{d-1}}{\lambda_2^{d-1}}|c|^d\,,
\end{equation}
for every $\lambda_2\in(0,1)$.

\bigbreak 

If $\lambda\in(0,1)$, we introduce a further parameter $\lambda_1\in(0,\lambda)$ so that
\[
\frac{\delta}{\e^{\lambda_1}}\to0 \hbox{ as } \e\to0.
\]

\smallbreak

\noindent Our construction relies on the definition of several concentric annuli. To this end, let 
\[
T:=\hbox{max}\,\{t\in \mathbb{N} : \e^{\lambda_1}2^t \leq R_{\Omega}\}=\lfloor\frac{\lambda_1 |\log \e|+\log R_{\Omega}}{\log 2}\rfloor
\]
and assume in particular that $T$ is larger than $4$ as $\e$ is small enough. Then pick a natural number $M\in(2,T)$ and define annuli centered in $z_{\e}$ having radii $\e^{\lambda_1}2^{kM}$, with $k=0,1,...,\lfloor\frac{T}{M}\rfloor+1$.

We have $\Omega\subseteq B(z_{\e},\e^{\lambda_1}2^{(\lfloor T/M\rfloor+1)M})$; hence, for every $u\in W^{1,d}_0(\Omega)$ such that $u=1$ on $B(z_{\e},\e)$, it holds that
\[
\begin{split}
F_{\e}(u,\Omega) & = F_{\e}(u,B(z_{\e},\e^{\lambda_1}2^{(\lfloor T/M\rfloor+1)M})) \\
& = F_{\e}(u,B(z_{\e},\e^{\lambda_2})) +  F_{\e}(u,B(z_{\e},\e^{\lambda_1})\setminus B(z_{\e},\e^{\lambda_2})) \\
& + \sum_{k=1}^{\lfloor T/M\rfloor+1}F_{\e}(u,B(z_{\e},\e^{\lambda_1}2^{kM})\setminus B(z_{\e},\e^{\lambda_1}2^{(k-1)M)}).
\end{split}
\]
In the last equality, we carefully separated three summands in order to treat each of them in accordance with the different exponential scales described by the parameters $\lambda_1, \lambda_2$. 

\noindent Apply Lemma \ref{modlemma} to the first summand  with 
\[
f(x,\xi)=f\left(\frac{x}{\delta},\xi\right),\, \eta=\e,\, R=\e^{\lambda_2},\, N\in \mathbb{N}\cap \left(1, \left \lfloor\frac{(1-\lambda_2)|\log\e|}{\log2} \right \rfloor \right) \hbox{ and } r=\e.
\]
Apply Lemma \ref{modlemma} to the second summand with 
\[
f(x,\xi)=f\left(\frac{x}{\delta},\xi\right),\, \eta=\e^{\lambda_2},\, R=\e^{\lambda_1},\, N\in\mathbb{N}\cap \left(1, \left \lfloor\frac{(\lambda_2-\lambda_1)|\log\e|}{\log2} \right \rfloor \right) \hbox{ and }  r=\e^{\lambda_2}.
\]
Apply Lemma \ref{modlemma} to the terms of the third summand for $k=1,...,\lfloor T/M\rfloor$ with 
\[
f(x,\xi)=f\left(\frac{x}{\delta},\xi\right),\, \eta=\e^{\lambda_1},\, R=\e^{\lambda_1}2^{kM},\, N\in\mathbb{N}\cap (1,kM) \hbox{ and }  r=\e^{\lambda_1}2^{(k-1)M}.
\]
Set for simplicity of notation
\[
S':= \left \lfloor\frac{(1-\lambda_2)|\log\e|}{\log2} \right \rfloor \qquad \hbox{and} \qquad S'':= \left \lfloor\frac{(\lambda_2-\lambda_1)|\log\e|}{\log2} \right \rfloor.
\]
Since $S', S''$ and $M$ will get arbitrarily large, we may assume we fix the same $N$ in each of the above applications of the lemma.

We obtain functions $v^{-1}\in W^{1,d}(B(z_{\e},\e^{\lambda_2}))$, $v^0\in W^{1,d}(B(z_{\e},\e^{\lambda_1})\setminus \overline{B}(z_{\e},\e^{\lambda_2}))$ and $v^k\in W^{1,d}(B(z,\e^{\lambda_1}2^{kM})\setminus \overline{B}(z,\e^{\lambda_1}2^{(k-1)M})),\, k=1,...,\lfloor T/M\rfloor$ with the properties stated in Lemma \ref{modlemma}. We put
\[
v:=
\begin{cases}
v^{-1} & \hbox{ on } B(z_{\e},\e^{\lambda_2}) \\
v^0 & \hbox{ on } B(z_{\e},\e^{\lambda_1})\setminus B(z_{\e},\e^{\lambda_2}) \\
v^k & \hbox{ on } B(z_{\e},\e^{\lambda_1}2^{kM})\setminus B(z_{\e},\e^{\lambda_1}2^{(k-1)M}),\, k=1,...,\lfloor T/M\rfloor \\
u & \hbox{ otherwise},
\end{cases}
\]
and note that $v\in W^{1,d}_0(B(z_{\e},{\e^{\lambda_1}2^{(\lfloor T/M\rfloor+1)M}}))$ since the modifications provided by the lemma occur far from the boundary of each annulus; moreover it holds
\[
\biggl(1+\frac{C}{N-1}\biggr)F_{\e}(u,B(z_{\e},\e^{\lambda_1}2^{(\lfloor T/M\rfloor+1)M})) \geq F_{\e}(v,B(z_{\e},\e^{\lambda_1}2^{(\lfloor T/M\rfloor+1)M})).
\]

In order to highlight that $v$ is constant of value $c_k$ on spheres centered in $z_{\e}$ of radii of the form $\e2^{S'-j_{-1}}$, $\e^{\lambda_2}2^{S''-j_0}$ and $\e^{\lambda_1}2^{kM-j_k}$, where $j_k\in\{1,...,N-1\}$ for $k=-1,0,1,...,\lfloor T/M\rfloor$, we write
\begin{equation}\label{below1lambdain}
\begin{split}
F_{\e}(v,B(z_{\e},\e^{\lambda_1}2^{(\lfloor T/M\rfloor+1)M})) & = F_{\e}(v,B(z_{\e},\e2^{S'-j_{-1}})) \\
& + F_{\e}(v,B(z_{\e},\e^{\lambda_2}2^{S''-j_0}) \setminus B(z_{\e},\e2^{S'-j_{-1}})) \\
& + F_{\e}(v, B(z_\e,\e^{\lambda_1}2^{M-j_1})\setminus B(z_{\e},\e^{\lambda_2}2^{S''-j_0})) \\
& + \sum_{k=2}^{\lfloor T/M\rfloor}F_{\e}(v,B(z_\e,\e^{\lambda_1}2^{kM-j_k})\setminus B(z_\e,\e^{\lambda_1}2^{(k-1)M-j_{k-1}})) \\
& + F_{\e}(v,B(z_{\e},\e^{\lambda_1}2^{(\lfloor T/M\rfloor+1)M})\setminus B(z_{\e},\e^{\lambda_1}2^{\lfloor T/M\rfloor M-j_{\lfloor T/M\rfloor}})).
\end{split}
\end{equation}
Then we define functions $w^k,\, k=-1,0,1,...\lfloor T/M \rfloor+1$ as follows: $w^{-1}\in W^{1,d}(B(z_{\e},\e^{\lambda_2}))$ is defined as
\[
w^{-1}:=
\begin{cases}
v & \hbox{ on } B(z_{\e},\e2^{S'-j_{-1}})\\
c_{-1} & \hbox{ otherwise},
\end{cases}
\]
so that 
\[
F_{\e}(w^{-1},B(z_{\e},\e^{\lambda_2})) = F_{\e}(v,B(z_{\e},\e2^{S'-j_{-1}})).
\]
Similarly, set
\[
w^{0}:=
\begin{cases}
c_{-1} & \hbox{ on } B(z_{\e}, \e2^{S'-j_{-1}}))\setminus \overline{B}(z_{\e}, \e2^{S'-N})) \\
v & \hbox{ on } B(z_{\e},\e^{\lambda_2}2^{S''-j_0}) \setminus B(z_{\e},\e2^{S'-j_{-1}})) \\
c_{0} & \hbox{ on } B(z_{\e}, \e^{\lambda_1}) \setminus B(z_{\e}, \e^{\lambda_2}2^{S''-j_0}),
\end{cases}
\]
so that 
\[
F_{\e}(w^{0},B(z_{\e},\e^{\lambda_1})\setminus B(z_{\e}, \e2^{S'-N}))) = F_{\e}(v,B(z_{\e},\e^{\lambda_2}2^{S''-j_0}) \setminus B(z_{\e},\e2^{S'-j_{-1}}))
\]
and
\[
w^{1}:=
\begin{cases}
c_{0} & \hbox{ on } B(z_{\e}, \e2^{S''-j_0})\setminus \overline{B}(z_{\e}, \e2^{S''-N})) \\
v & \hbox{ on } B(z_{\e},\e^{\lambda_1}2^{M-j_1}) \setminus B(z_{\e},\e2^{S''-j_0})) \\
c_{1} & \hbox{ on } B(z_{\e}, \e^{\lambda_1}2^{M})) \setminus B(z_{\e}, \e^{\lambda_1}2^{M-j_1}),
\end{cases}
\]
so that
\[
F_{\e}(w^{1},B(z_{\e},\e^{\lambda_1}2^M)\setminus B(z_{\e}, \e2^{S''-N}))) = F_{\e}(v,B(z_{\e},\e^{\lambda_1}2^{M-j_1}) \setminus B(z_{\e},\e2^{S''-j_{0}})).
\]
For $k=2,...,\lfloor T/M \rfloor+1$, we define annuli
\[
A_{M,k}^N:=B(z_{\e}, \e^{\lambda_1}2^{kM}) \setminus \overline{B}(z_{\e}, \e^{\lambda_1}2^{(k-1)M-N}).
\]
For $k=2,..., \lfloor T/M\rfloor$, we define functions $w^k\in W^{1,d}(A_{M,k}^N)$ as
\[
w^{k}:=
\begin{cases}
c_{k-1} & \hbox{ on } B(z_{\e},\e^{\lambda_1}2^{(k-1)M-j_{k-1}})\setminus \overline{B}(z_{\e},\e^{\lambda_1}2^{(k-1)M-N}) \\
v & \hbox{ on } B(z_{\e},\e^{\lambda_1}2^{kM-j_k})\setminus B(z,\e^{\lambda_1}2^{(k-1)M-j_{k-1}}) \\
c_{k} & \hbox{ on } B(z,\e^{\lambda_1}2^{kM}) \setminus B(z,\e^{\lambda_1}2^{kM-j_{k}}),
\end{cases}
\]
and for $k=\lfloor T/M \rfloor+1$,
\[
w^{\lfloor T/M \rfloor+1} :=
\begin{cases}
c_{\lfloor T/M \rfloor} & \hbox{ on } B(z_{\e},\e^{\lambda_1}2^{(k-1)M-j_{k-1}})\setminus \overline{B}(z_{\e},\e^{\lambda_1}2^{(k-1)M-N}) \\
v & \hbox{ otherwise},
\end{cases}
\]
so that
\[
F_{\e}(w^k, A_{M,k}^N) = F_{\e}(v,B(z,\e^{\lambda_1}2^{kM-j_k})\setminus B(z,\e^{\lambda_1}2^{(k-1)M-j_{k-1}}))
\]
for all $k=2,...,\lfloor T/M \rfloor+1$.

\smallbreak

If we set $A_{M,-1}^N:= B(z_{\e},\e^{\lambda_2})$, $A_{M,0}^N:=B(z_{\e},\e^{\lambda_1}) \setminus \overline{B}(z_{\e},\e2^{S'-N})$ and $A_{M,1}^N:=B(z_{\e},\e^{\lambda_1}2^M)\setminus \overline{B}(z_{\e}, \e2^{S''-N})$, then we can rewrite \eqref{below1lambdain} simply as
\[
F_{\e}(v,B(z_\e,\e^{\lambda_1}2^{(\lfloor T/M\rfloor+1)M})) =  \sum_{k=-1}^{\lfloor T/M\rfloor+1} F_{\e}(w^k,A_{M,k}^N).
\]
Once more, we take advantage of the fact that the functions $w^{-1},...,w^{\lfloor T/M\rfloor+1}$ attain constant value on the components of their annuli of definition. Also, exploiting (P$2$) and suitable affine transformation (as in the case $\lambda=0$), we get 
\[
\biggl(1+\frac{C}{N-1}\biggr)F_{\e}(u,\Omega) \geq 
\]
\begin{align}
& \geq \min\{F_{\e}(v,B(z_{\e},\e^{\lambda_2})) : v\in W^{1,d}_0(B(z_{\e},\e^{\lambda_2})), v=1 \hbox{ on } B(z_{\e},\e)\}|1-c_{-1}|^d \label{in01-1} \\
& + \min\{F_{\e}(v,B(z_{\e},\e^{\lambda_1})  \setminus B(z_{\e},\e2^{S'-N})) : v\in W^{1,d}(B(z_{\e},\e^{\lambda_1})\setminus \overline{B}(z_\e,\e2^{S'-N})), \notag \\ 
& \qquad \qquad \qquad \qquad v=1 \hbox{ on } \partial B(z_{\e},\e2^{S'-N}), v=0 \hbox{ on } \partial B(z_{\e},\e^{\lambda_1})\}|c_{-1}-c_0|^d \label{in01-2} \\
& + \min\{F_{\e}(v,B(z_{\e},\e^{\lambda_1}2^M) \setminus B(z_{\e},\e2^{S''-N})) : v\in W^{1,d}(B(z_{\e},\e^{\lambda_1}2^M) \setminus \overline{B}(z_{\e},\e2^{S''-N})), \notag \\ 
& \qquad \qquad \qquad \qquad v=1 \hbox{ on } \partial B(z_{\e},\e2^{S''-N}), v=0 \hbox{ on } \partial B(z_{\e},\e^{\lambda_1}2^M)\}|c_{0}-c_1|^d \label{in01-3} \\
& + \sum_{k=2}^{\lfloor T/M\rfloor+1} \min\{F_{\e}(v, A_{M,k}^N) : v\in W^{1,d}(A_{M,k}^N), \,  v=1 \hbox{ on } \partial B(z_{\e},\e^{\lambda_1}2^{(k-1)M-N}), \notag \\
& \qquad \qquad \qquad \qquad v=0 \hbox{ on } \partial B(z_{\e},\e^{\lambda_1}2^{kM})\}|c_{k-1}-c_{k}|^d\,, \label{in01-4}
\end{align}

where we put $c_{\lfloor \frac{T}{M} \rfloor+1}:=0$.

\bigbreak

Since $\lambda_2<\lambda$, the minimum in \eqref{in01-1} is estimated as for \eqref{below2lambda0} in the case $\lambda=0$, thus it is greater than or equal to
\begin{equation}\label{below2lambdain}
\frac{\Phi(z)+o_{\e}(1)}{(1-\lambda_2)^{d-1}|\log\e|^{d-1}}|1-c_{-1}|^d.
\end{equation}

The bounds for the \eqref{in01-2} and \eqref{in01-3} follow again by the growth condition from below in (P$3$), in particular, recalling how we defined $S'$ and $S''$, we have
\begin{equation}\label{below3lambdain}
\begin{split}
\alpha\hbox{Cap}_d(B(z_{\e},\e2^{S'-N}),B(z_{\e},\e^{\lambda_1})) & \geq \frac{\alpha\sigma_{d-1}}{[(1-\lambda_1)|\log\e|-(S'-N)\log2]^{d-1}} \\
& \geq \frac{\alpha\sigma_{d-1}}{[(\lambda_2-\lambda_1)|\log\e|+(N+1)\log2]^{d-1}}
\end{split}
\end{equation}
while
\begin{equation}\label{below4lambdain}
\begin{split}
\alpha\hbox{Cap}_d(B(z_{\e},\e2^{S''-N}),B(z_{\e},\e^{\lambda_1}2^M)) & \geq \frac{\alpha\sigma_{d-1}}{[M\log2+(1-\lambda_1)|\log\e|-(S''-N)\log2]^{d-1}} \\
& \geq \frac{\alpha\sigma_{d-1}}{[M\log2+(1-\lambda_2)|\log\e|+(N+1)\log2]^{d-1}}\,.
\end{split}
\end{equation}

Concerning the summands in \eqref{in01-4}, we proceed fixing $k=2,...,\lfloor T/M\rfloor+1$ and applying $v(x)\mapsto v(z_{\e}+x\e^{\lambda_1}2^{(k-1)M-N})$, so that each term equals
\[
\begin{split}
\min\Bigl\{\int_{B(0,2^{M+N})\setminus B(0,1)} & f\Bigl(\frac{x}{\delta}\e^{\lambda_1}2^{(k-1)M-N}+\frac{z_{\e}}{\delta},\nabla v(x)\Bigr)\,dx : \\
& v\in W^{1,d}_0(B(0,2^{M+N})), v=1 \hbox{ on } B(0,1)\Bigr\}|c_{k-1}-c_{k}|^d\,.
\end{split}
\]
By $\lambda_1<\lambda$ it follows that 
\[
\frac{\delta}{\e^{\lambda_1}2^{(k-1)M-N}}\to0 \hbox{ as } \e\to0;
\]
hence, we can apply Theorem \ref{homthmtranslation} with 
\[
A=B(0,2^{M+N})\setminus \overline{B}(0,1)\,, \qquad \eta=\frac{\delta}{\e^{ \lambda_1}2^{(k-1)M-N}}\,, \qquad \tau_{\eta}=\frac{z_{\e}}{\delta} 
\]
and $\phi$ any function in $W^{1,d}(B(0,2^{M+N})\setminus \overline{B}(0,1))$
such that $\phi=1$ on $\partial B(0,1)$ and $\phi=0$ on $\partial B(0,2^{M+N})$.
We get that each of the above minima equals
\begin{equation}\label{below5lambdain}
\begin{split}
\biggl[\min\Bigl\{\int_{B(0,2^{M+N})\setminus B(0,1)} & f_{\rm hom}(\nabla v(x))\,dx:\, v\in W^{1,d}(B(0,2^{M+N})), \\ 
& v=1 \hbox{ on } B(0,1), v=0 \hbox{ on } \partial B(0,2^{M+N})\Bigr\} + o_{\e}(1)\biggr]|c_{k-1}-c_{k}|^d, 
\end{split}
\end{equation}
where $f_{\rm hom}$ is the $d$-homogeneous function given by \eqref{fhom}, which does not depend on $k$.

\noindent Recalling the definition of the constant $C_{\rm hom}$ given in \eqref{Chom}, \eqref{below5lambdain} turns into
\[
\biggl[\frac{C_{\rm hom}+o_{M}(1)}{((M+N)\log2)^{d-1}}+o_{\e}(1)\biggr]|c_{k-1}-c_{k}|^d.
\]
We sum over $k$ and use the convexity of $x\mapsto |x|^d$, the fact that $\sum_{k=2}^{\lfloor T/M \rfloor+1} (c_{k-1}-c_k)=c_1$ and that $T\leq \frac{\lambda_1|\log\e|+\log R_{\Omega}}{\log2}$; we obtain
\[
\sum_{k=2}^{\lfloor T/M\rfloor+1} |c_{k-1}-c_{k}|^d 
\geq \frac{ (M\log2)^{d-1}}{(\lambda_1|\log\e|+\log R_{\Omega}+M\log2)^{d-1}}|c_1|^d,
\]
and in turn
\begin{multline}\label{below6lambdain}
\sum_{k=2}^{\lfloor T/M\rfloor+1} \biggl[\frac{C_{\rm hom}+o_{M}(1)}{((M+N)\log2)^{d-1}}+o_{\e}(1)\biggr]|c_{k-1}-c_{k}|^d \\ \geq \biggl[\frac{C_{\rm hom}+o_{M}(1)}{((M+N)\log2)^{d-1}}+o_{\e}(1)\biggr]\frac{ (M\log2)^{d-1}}{(\lambda_1|\log\e|+\log R_{\Omega}+M\log2)^{d-1}}|c_1|^d\,.
\end{multline}

Gathering \eqref{below2lambdain}, \eqref{below3lambdain}, \eqref{below4lambdain} and \eqref{below6lambdain}, and multiplying by $|\log\e|^{d-1}$, we get
\begin{multline*}
\biggl(1+\frac{C}{N-1}\biggr)|\log\e|^{d-1}F_{\e}(u,\Omega) 
\geq \frac{\Phi(z)+o_{\e}(1)}{(1-\lambda_2)^{d-1}}|1-c_{-1}|^d \\
+ \frac{\alpha\sigma_{d-1}|\log\e|^{d-1}}{[(\lambda_2-\lambda_1)|\log\e|+(N+1)\log2]^{d-1}}|c_{-1}-c_0|^d \\
+ \frac{\alpha\sigma_{d-1}|\log\e|^{d-1}}{[M\log2+(1-\lambda_2)|\log\e|+(N+1)\log2]^{d-1}}|c_{0}-c_1|^d \\
+ \biggl[\frac{C_{hom}+o_{M}(1)}{((M+N)\log2)^{d-1}}+o_{\e}(1)\biggr]\frac{ |\log\e|^{d-1}(M\log2)^{d-1}}{(\lambda_1|\log\e|+\log R_{\Omega}+M\log2)^{d-1}}|c_1|^d\,.
\end{multline*}

\medbreak

Arguing as before, we stress that $c_{-1}, c_0$ and $c_1$ depend on $\e$ and can be picked inside the interval $[0,1]$. This lead us to assume that they all converge to some finite limits, say $c_{-1}, c_0, c_1$, respectively. Moreover, such limits have to coincide; if not we would get a contradiction letting $\lambda_1,\lambda_2\to\lambda$ or $\lambda_2\to1$. 

\noindent Eventually, the following estimate holds true:
\[
\begin{split}
\biggl(1+\frac{C}{N-1}\biggr)\liminf_{\e\to0}|\log\e|^{d-1}F_{\e}(u,\Omega) & 
\geq \frac{\Phi(z)}{(1-\lambda_2)^{d-1}}|1-c|^d \\
& + \biggl[\frac{C_{\rm hom}+o_{M}(1)}{((M+N)\log2)^{d-1}}\biggr]\frac{ (M\log2)^{d-1}}{\lambda_1^{d-1}}|c|^d
\end{split}
\]
and letting $M\to+\infty, N\to+\infty$, by the arbitrariness of $u$ we achieve 
\begin{equation}\label{importantliminf2}
\liminf_{\e\to0}|\log\e|^{d-1}\mu_{\e,\delta} \geq \frac{\Phi(z)}{(1-\lambda_2)^{d-1}}|1-c|^d + \frac{C_{\rm hom}}{\lambda_1^{d-1}}|c|^d\,.
\end{equation}

\bigbreak

If $\lambda=1$, keeping the notation introduced throughout the proof, define annuli centered in $z_{\e}$ having radii $\e^{\lambda_1}2^{kM}$, with $k=1,...,\lfloor\frac{T}{M}\rfloor+1$.

\noindent For every function $u\in W^{1,d}_0(\Omega), u=1$ on $B(z_{\e},\e)$, we have
\begin{equation}\label{mainpropfixedhom}
\begin{split}
F_{\e}(u,\Omega) 
& = F_{\e}(u,B(z_{\e},\e^{\lambda_1}2^{(\lfloor T/M\rfloor+1)M})) \\
& =  F_{\e}(u,B(z_{\e},\e^{\lambda_1})) \\
& + \sum_{k=1}^{\lfloor T/M\rfloor+1}F_{\e}(u,B(z_{\e},\e^{\lambda_1}2^{kM})\setminus B(z_{\e},\e^{\lambda_1}2^{(k-1)M})).
\end{split}
\end{equation}
Apply Lemma \ref{modlemma} to the terms of the second summand for $k=1,...,\lfloor T/M\rfloor$ with 
\[
f(x,\xi)=f\left(\frac{x}{\delta},\xi\right),\, \eta=\e^{\lambda_1},\, R=\e^{\lambda_1}2^{kM},\, N\in\mathbb{N}\cap(1,kM) \hbox{ and } r=\e^{\lambda_1}2^{(k-1)M}.
\]
Arguing as in the previous instances, with $\lambda\in[0,1)$, we get
\[
\biggl(1+\frac{C}{N-1}\biggr)F_{\e}(u,\Omega) \geq
\]
\begin{align}
& \geq \min\{F_{\e}(v,B(z_{\e},\e^{\lambda_1})) : v\in W^{1,d}(B(z_{\e},\e^{\lambda_1})), v=1 \hbox{ on } B(z_{\e},\e), \notag \\
& \qquad \qquad \qquad \qquad \qquad \qquad \qquad v=0 \hbox{ on } \partial B(z_{\e},\e^{\lambda_1})\}|1-c_0|^d \label{in1-1} \\ 
& + \sum_{k=1}^{\lfloor T/M\rfloor+1} \min\{F_{\e}(v, A_{M,k}^N) : v\in W^{1,d}(A_{M,k}^N), v=1 \hbox{ on } \partial B(z_{\e},\e^{\lambda_1}2^{(k-1)M-N}), \notag \\
& \qquad \qquad \qquad \qquad \qquad \qquad \qquad  v=0 \hbox{ on } \partial B(z_{\e},\e^{\lambda_1}2^{kM})\}|c_{k-1}-c_{k}|^d\,, \label{in1-2}
\end{align}
where $c_{\lfloor \frac{T}{M} \rfloor+1}:=0$.

Making use of (P$3$), \eqref{in1-1} is bounded from below by
\[
\alpha\hbox{Cap}_d(B(z_{\e},\e),B(z_{\e},\e^{\lambda_1}))|1-c_0|^d = \frac{\alpha\sigma_{d-1}}{[(1-\lambda_1)|\log\e|]^{d-1}}|1-c_0|^d\,,
\]
while \eqref{in1-2} can be estimated as in \eqref{below6lambdain} since $\delta/\e^{\lambda_1}\to0$.

At the end, we get the inequality
\begin{multline*}
\biggl(1+\frac{C}{N-1}\biggr)|\log\e|^{d-1}F_{\e}(u,\Omega) \geq \frac{\alpha\sigma_{d-1}}{(1-\lambda_1)^{d-1}}|1-c_0|^d \\
+ \biggl[\frac{C_{\rm hom}+o_{M}(1)}{((M+N)\log2)^{d-1}}+o_{\e}(1)\biggr]\frac{ |\log\e|^{d-1}(M\log2)^{d-1}}{(\lambda_1|\log\e|+\log R_{\Omega}+M\log2)^{d-1}}|c_0|^d\,.
\end{multline*}

Recall that we may assume that $c_0=c_0(\e)$ converges to a finite value $c$, hence we let $\e\to0, M\to+\infty$ and $N\to+\infty$, then we take advantage of the arbitrariness of $u$, to obtain
\begin{equation}\label{importantliminf3}
\liminf_{\e\to0}|\log\e|^{d-1}\mu_{\e,\delta} \geq \frac{\alpha\sigma_{d-1}}{(1-\lambda_1)^{d-1}}|1-c|^d + \frac{C_{\rm hom}}{\lambda_1^{d-1}}|c|^d\,.
\end{equation}

\bigbreak

Once we gather \eqref{importantliminf1}, \eqref{importantliminf2}, \eqref{importantliminf3}, we have 
\[
\liminf_{\e\to0}|\log\e|^{d-1}\mu_{\e,\delta} \geq
\begin{cases}
\frac{\Phi(z)}{(1-\lambda_2)^{d-1}}|1-c|^d + \frac{\alpha\sigma_{d-1}}{\lambda_2^{d-1}}|c|^d & \hbox{ if } \lambda=0, \\[5pt]
\frac{\Phi(z)}{(1-\lambda_2)^{d-1}}|1-c|^d + \frac{C_{\rm hom}}{\lambda_1^{d-1}}|c|^d & \hbox{ if } \lambda\in(0,1), \\[5pt]
\frac{\alpha\sigma_{d-1}}{(1-\lambda_1)^{d-1}}|1-c|^d + \frac{C_{\rm hom}}{\lambda_1^{d-1}}|c|^d & \hbox{ if } \lambda=1
\end{cases}
\]
for every $\lambda_1\in(0,\lambda)$ and $\lambda_2\in(\lambda,1)$.

\noindent These expressions can be estimated by the same argument concerning the minimization of the function $a|1-x|^d+b|x|^d$ with $a,b>0$. Indeed, the minimum is attained at 
\[
x=\biggl[\Bigl(\frac{b}{a}\Bigr)^{\frac{1}{d-1}}+1\biggr]^{-1}
\]
with minimum value 
\[
b\biggl[\Bigl(\frac{b}{a}\Bigr)^{\frac{1}{d-1}}+1\biggr]^{1-d}.
\]

In \eqref{importantliminf1}, we set
\[
a=\frac{\Phi(z)}{(1-\lambda_2)^{d-1}} \qquad \hbox{and} \qquad b=\frac{\alpha\sigma_{d-1}}{\lambda_2^{d-1}}\,,
\]
to achieve
\[
\begin{split}
\liminf_{\e\to0}|\log\e|^{d-1}\mu_{\e,\delta} & \geq 
\frac{\alpha\sigma_{d-1}}{\lambda_2^{d-1}}\left[\left(\frac{\alpha\sigma_{d-1}/\lambda_2^{d-1}}{\Phi(z)/(1-\lambda_2)^{d-1}}\right)^{\frac{1}{d-1}}+1 \right]^{1-d} \\
& = \Phi(z)\alpha\sigma_{d-1}\left[(1-\lambda_2)(\alpha\sigma_{d-1})^{\frac{1}{d-1}}+\lambda_2\Phi(z)^{\frac{1}{d-1}}\right]^{1-d}.
\end{split}
\]
We conclude passing to the limit as $\lambda_2\to0$.

In \eqref{importantliminf2}, put
\begin{equation} \label{1dfunction}
    a=\frac{\Phi(z)}{(1-\lambda_2)^{d-1}} \qquad \hbox{and} \qquad b=\frac{C_{\rm hom}}{\lambda_1^{d-1}}\,,
\end{equation}
and let $\lambda_1,\lambda_2\to\lambda$ getting
\[
\begin{split}
\liminf_{\e\to0}|\log\e|^{d-1}\mu_{\e,\delta} & \geq \frac{C_{\rm hom}}{\lambda^{d-1}}\left[\left(\frac{C_{\rm hom}/\lambda^{d-1}}{\Phi(z)/(1-\lambda)^{d-1}}\right)^{\frac{1}{d-1}}+1 \right]^{1-d}  \\
& = \Phi(z)C_{\rm hom}\Bigl[(1-\lambda)C_{\rm hom}^{\frac{1}{d-1}}+\lambda\Phi(z)^{\frac{1}{d-1}}\Bigr]^{1-d}.
\end{split}
\]

Finally, in \eqref{importantliminf3} let 
\[
a = \frac{\alpha\sigma_{d-1}}{(1-\lambda_1)^{d-1}} \qquad \hbox{and} \qquad b = \frac{C_{\rm hom}}{\lambda_1^{d-1}}\,,
\]
to have
\[
\begin{split}
\liminf_{\e\to0}|\log\e|^{d-1}\mu_{\e,\delta} & \geq \frac{C_{\rm hom}}{\lambda_1^{d-1}}\left[\left(\frac{C_{\rm hom}/\lambda_1}{\alpha\sigma_{d-1}/(1-\lambda_1)}\right)^{\frac{1}{d-1}}+1\right]^{1-d} \\
& = \alpha\sigma_{d-1} C_{\rm hom}\Bigl[(1-\lambda_1)C_{\rm hom}^{\frac{1}{d-1}}+\lambda_1(\alpha\sigma_{d-1})^{\frac{1}{d-1}}\Bigr]^{1-d}.
\end{split}
\]
Then, conclude letting $\lambda_1\to1$.

\subsection{Construction of optimal sequences}

To finish the proof we define minimizing sequences providing the bound from above using suitable capacitary profiles. 

\medbreak

If $\lambda=0$, take $\lambda_2\in(\lambda,1)$ and let $v_{\e}^0$ be a solution of the minimum problem
\[
\min\Bigl\{\int_{B(0,\e^{\lambda_2-1})} f(z, \nabla u(x))\,dx : u\in W^{1,d}_0(B(0,\e^{\lambda_2-1})), u=1 \hbox{ on } B(0,1)\Bigr\}\,.
\] 

For $\e\ll1$, the function $u_{\e}^0(x):=v_{\e}^0\left(\frac{x-z_{\e}}{\e}\right)$ belongs to $W^{1,d}_0(\Omega)$ and it is admissible for the minimum problem defining \eqref{minpbm2}, thus $\mu_{\e,\delta} \leq F_{\e}(u_{\e}^0,\Omega)$ and
by change of variables and homogeneity of $f(x,\cdot)$, it holds
\[
\begin{split}
F_{\e}(u_{\e}^0,\Omega) & = F_{\e}(u_{\e}^0,B(z_\e, \e^{\lambda_2})) = \int_{B(z_{\e},\e^{\lambda_2})} f\Bigl(\frac{x}{\delta},\nabla u_{\e}^0(x)\Bigr)\,dx \\
& = \int_{B(0,\e^{\lambda_2-1})} f\Bigl(\frac{x}{\delta}\e+\frac{z_{\e}}{\delta},\nabla v_{\e}^0(x)\Bigr)\,dx = \int_{B(0,\e^{\lambda_2-1})} f\Bigl(\frac{x}{\delta}\e + z,\nabla v_{\e}^0(x)\Bigr)\,dx.
\end{split}
\]
Note that $\e^{\lambda_2}/\delta\to0$, hence, for every $x\in B(0,\e^{\lambda_2-1})$, it holds $|x|\frac{\e}{\delta}\to0$ as $\e\to0$. In light of this, given any $\nu>0$, by \eqref{assumption0} we deduce that
\[
f(x,\xi) \leq f(z,\xi)+\nu |\xi|^d \text{ for every } \xi\in\mathbb{R}^d \text{ and for every } x\in B(0,\e^{\lambda_2-1}),
\]
as $\e$ is sufficiently small. As a consequence 
\[
\int_{B(0,\e^{\lambda_2-1})} f\Bigl(\frac{x}{\delta}\e+z,\nabla v_{\e}^0(x)\Bigr)\,dx \leq \int_{B(0,\e^{\lambda_2-1})} f(z,\nabla v_{\e}^0(x))\,dx + \nu\int_{B(0,\e^{\lambda_2-1})} |\nabla v_{\e}^0(x)|^d\,dx 
\]
which, by the growth condition, is bounded above by
\[
\left(1+\frac{\nu}{\alpha}\right)\int_{B(0,\e^{\lambda_2-1})} f(z,\nabla v_{\e}^0(x))\,dx \,.
\]
In light of the fact that $\e^{\lambda_2-1}\to\infty$ as $\e\to0$, we apply Lemma \ref{limlemma} to deduce
\begin{equation}\label{above1}
F_{\e}(u_{\e}^0,\Omega) = F_{\e}(u_{\e}^0,B(z_\e, \e^{\lambda_2})) \leq \left(1+\frac{\nu}{\alpha}\right)\frac{\Phi(z)+o_{\e}(1)}{(1-\lambda_2)^{d-1}|\log\e|^{d-1}}\,.
\end{equation}
Thus, we conclude by the arbitrariness of $\nu>0$ and $\lambda_2\in(0,1)$, that 
\[
\limsup_{\e\to0}|\log\e|^{d-1}\mu_{\e,\delta} \leq \inf_{\lambda_2\in(0,1)}\frac{\Phi(z)}{(1-\lambda_2)^{d-1}}=\Phi(z).
\]

\bigbreak

If $\lambda\in(0,1)$, introduce a further parameter $\lambda_1\in(0,\lambda)$, put
\[
T:=\hbox{max}\{t\in\mathbb{N} : \e^{\lambda_1}2^t\leq \text{dist}(z_{\e},\partial\Omega)\} = \left \lfloor \frac{\lambda_1|\log\e|+\log \text{dist}(z_{\e},\partial\Omega)}{\log2} \right \rfloor
\]
and take $M\in\mathbb{N}\cap (0,T)$. Since the family of points $\{z_{\e},\, \e>0\}$ is contained in a ball, say $B$, whose closure lays inside $\Omega$, we have that $ \text{dist}(z_{\e},\partial \Omega) \geq \text{dist}(\partial B,\partial \Omega)>0$ so that $T$ is well defined and can be assumed to be greater than $2$ for every $\e$.

\smallbreak

Let $v_{\eta}$ be a solution of the minimum problem
\[
m_{\eta}:=\min \Bigl\{\int_{B(0,2^{M})} f\Bigl(\frac{x}{\eta}+ \tau_{\eta}, \nabla u(x)\Bigr)\,dx : \\ u\in W^{1,d}_0(B(0,2^{M})), u=1 \hbox{ on } B(0,1)\Bigr\}\,
\] 
and set
\[
m_0:=\min \Bigl\{\int_{B(0,2^{M})} f_{\rm hom}(\nabla u(x))\,dx : \\ u\in W^{1,d}_0(B(0,2^{M})), u=1 \hbox{ on } B(0,1)\Bigr\}.
\]
By Theorem \ref{homthmtranslation}, there exists an increasing non negative function $\omega$ such that
\[
|m_{\eta}-m_0| \leq \omega(\eta)\,\, \hbox{ and }\,\, \omega(\eta)\to0 \hbox{ as } \eta\to0;
\]
thus, for $k=1,...,\lfloor T/M \rfloor$, define $u_{\e}^k\in W^{1,d}(B(z_{\e}, \e^{\lambda_1}2^{kM})\setminus \overline{B}(z_{\e}, \e^{\lambda_1}2^{(k-1)M}))$ as 
\[
u_{\e}^k(x):=\frac{c}{\lfloor T/M\rfloor} v_{\eta}\left(\frac{x-z_{\e}}{\e^{\lambda_1}2^{(k-1)M}}\right) + \frac{\lfloor T/M\rfloor-k}{\lfloor T/M\rfloor}c\,
\]
for some constant $c$ to be properly selected.

If we set $\eta=\frac{\delta}{\e^{\lambda_1}2^{(k-1)M}}, \tau_{\eta}=\frac{z_{\e}}{\delta}$ and we apply a change of variables and the homogeneity of $f(x,\cdot)$, it holds
\begin{equation}\label{etad}
\begin{split}
F_{\e}(u_{\e}^k,B(z_{\e} &,\e^{\lambda_1}2^{kM}) \setminus B(z_{\e},\e^{\lambda_1}2^{(k-1)M})) \\
& = \int_{B(z_{\e},\e^{\lambda_1}2^{kM})\setminus B(z_{\e},\e^{\lambda_1}2^{(k-1)M})} f\left(\frac{x}{\delta}, \nabla u_{\e}^k(x)\right)\,dx \\
& = \left|\frac{c}{\lfloor T/M \rfloor}\right|^d \int_{B(0,2^M)} f\left(\frac{x}{\delta}\e^{\lambda_1}2^{k-1}+\frac{z_{\e}}{\delta}, \nabla v_{\eta}(x)\right)\,dx \\
& = \left|\frac{c}{\lfloor T/M \rfloor}\right|^d m_{\eta}
\\
& \leq \left|\frac{c}{\lfloor T/M \rfloor}\right|^d (m_0+\omega(\eta)) \\
& \leq \left|\frac{c}{\lfloor T/M \rfloor}\right|^d \left(m_0+\omega\left(\frac{\delta}{\e^{\lambda_1}}\right)\right).
\end{split}
\end{equation}

\noindent Then, considering the same $u_{\e}^0$ introduced in the case $\lambda=0$, set
\[
u_{\e}(x):=
\begin{cases}
(1-c)u_{\e}^0(x)+c & \hbox{ if } x\in B(z_{\e},\e^{\lambda_2}) \\
c & \hbox{ if } x\in B(z_{\e},\e^{\lambda_1})\setminus B(z_{\e},\e^{\lambda_2}) \\
u_{\e}^k(x) & \hbox{ if } x\in B(z_{\e},\e^{\lambda_1}2^{kM})\setminus B(z_{\e},\e^{\lambda_1}2^{(k-1)M})\,, k=1,...,\lfloor T/M \rfloor \\
0 & \hbox{ if } x\in \Omega\setminus B(z_{\e}, \e^{\lambda_1}2^{\lfloor T/M\rfloor M}).
\end{cases}
\]
Since the boundary conditions match, $u_{\e}\in W^{1,d}_0(\Omega)$ and $u_{\e}=1$ on $B(z_{\e},\e)$; therefore, it is an admissible function for the minimum problem.

We estimate $F_{\e}(u_{\e},B(z_{\e},\e^{\lambda_1}))$ and $F_{\e}(u_{\e},\Omega\setminus B(z_{\e},\e^{\lambda_1}))$ separately, neglecting those regions on which $u_{\e}$ is constant.

\smallbreak

By the same computation which led to \eqref{above1},
\begin{equation}\label{above2}
F_{\e}(u_{\e},B(z_{\e}, \e^{\lambda_1})) = F_{\e}(u_{\e}^0,B(z_{\e}, \e^{\lambda_2}))|1-c|^d \leq \frac{\Phi(z)+o_{\e}(1)}{(1-\lambda_2)^{d-1}|\log\e|^{d-1}}|1-c|^d\,,
\end{equation}
where $\nu$ has been neglected as it gets arbitrarily small as $\e\to0$.

\smallbreak

We focus on $F_{\e}(u_{\e},\Omega\setminus B(z_{\e}, \e^{\lambda_1}))$. By \eqref{etad}, it holds 
\begin{multline}\label{above3}
F_{\e}(u_{\e},\Omega\setminus B(z_{\e}, \e^{\lambda_1})) = \sum_{k=1}^{\lfloor T/M\rfloor} F_{\e}(u_{\e}^k,B(z_{\e},\e^{\lambda_1}2^{kM})\setminus B(z_{\e},\e^{\lambda_1}2^{(k-1)M})) \\
\leq \left|\frac{c}{\lfloor T/M \rfloor}\right|^d \sum_{k=1}^{\lfloor T/M\rfloor} \left(m_0+\omega\left(\frac{\delta}{\e^{\lambda_1}}\right)\right) 
= \frac{|c|^d}{\lfloor T/M \rfloor^{d-1}} \left(m_0+\omega\left(\frac{\delta}{\e^{\lambda_1}}\right)\right),
\end{multline}
but $f_{\rm hom}(0)=0$, while by the definition of the constant $C_{\rm hom}$, we have
\[
m_0=\frac{C_{\rm hom}+o_M(1)}{(M \log 2)^{d-1}}\,.
\]
Thus, we substitute in \eqref{above3} obtaining
\[
F_{\e}(u_{\e},\Omega\setminus B(z_{\e}, \e^{\lambda_1})) \leq \frac{|c|^d}{\lfloor T/M \rfloor^{d-1}} \left(\frac{C_{\rm hom}+o_M(1)}{(M \log 2)^{d-1}} + \omega\left(\frac{\delta}{\e^{\lambda_1}}\right)\right),
\]
and, as $T \geq \frac{\lambda_1|\log\e|+\log d - \log2}{\log2}$ with $d:=$dist$(\partial B,\partial\Omega)$, it holds
\begin{equation}\label{above4}
\begin{split}
F_{\e}(u_{\e},\Omega\setminus B(z_{\e}, \e^{\lambda_1})) \leq \frac{C_{\rm hom}+o_M(1) + (M\log2)^{d-1} \omega(\delta/\e^{\lambda_1})}{(\lambda_1|\log\e|+\log d-\log2)^{d-1}}|c|^d.
\end{split}    
\end{equation}

\smallbreak

We gather estimates \eqref{above2}, \eqref{above4} to conclude
\[
\begin{split}
|\log\e|^{d-1} \mu_{\e,\delta} & \leq \frac{\Phi(z)+o_{\e}(1)}{(1-\lambda_2)^{d-1}}|1-c|^d \\
& + \frac{|\log\e|^{d-1}[C_{\rm hom}+o_M(1) + (M\log2)^{d-1}\omega(\delta/\e^{\lambda_1})] }{(\lambda_1|\log\e|+\log d-\log2)^{d-1}}|c|^d\,.
\end{split}
\]

\noindent Since $\frac{\delta}{\e^{\lambda_1}}\to0$, let $\e\to0$ and then $M\to+\infty$ to deduce
\[
\limsup_{\e\to0}|\log\e|^{d-1}\mu_{\e,\delta} \leq \frac{\Phi(z)}{(1-\lambda_2)^{d-1}}|1-c|^d +  \frac{C_{\rm hom}}{\lambda_1^{d-1}}|c|^d\,;
\]
then, let $\lambda_1,\lambda_2\to\lambda$ so that
\[
\limsup_{\e\to0}|\log\e|^{d-1}\mu_{\e,\delta} \leq \frac{\Phi(z)}{(1-\lambda)^{d-1}}|1-c|^d +  \frac{C_{\rm hom}+o_M(1)}{\lambda^{d-1}}|c|^d.
\]

Finally, put $c:=\Bigl[\Bigl(\frac{b}{a}\Bigr)^{\frac{1}{d-1}}+1\Bigr]^{-1}$,
with $a=\Phi(z)/(1-\lambda)^{d-1},\, b=C_{\rm hom}/\lambda^{d-1}$. As we are exactly in the case discussed in  \eqref{1dfunction} with $\lambda=\lambda_1=\lambda_2$, the same computation holds, leading to
\[
\limsup_{\e\to0}|\log\e|^{d-1}\mu_{\e,\delta} \leq \Phi(z)C_{\rm hom}\Bigl[(1-\lambda)C_{\rm hom}^{\frac{1}{d-1}}+\lambda\Phi(z)^{\frac{1}{d-1}}\Bigr]^{1-d}.
\]

\bigbreak

If $\lambda=1$ we just set $c=1$ and
\[
u_{\e}(x):=
\begin{cases}
1 & \hbox{ if } x\in B(z_{\e},\e^{\lambda_1}) \\
u_{\e}^k(x) & \hbox{ if } x\in B(z_{\e},\e^{\lambda_1}2^{kM})\setminus B(z_{\e},\e^{\lambda_1}2^{(k-1)M}), k=1,...,\lfloor T/M\rfloor \\
0 & \hbox{ if } x\in\Omega\setminus B(z_{\e}, \e^{\lambda_1}2^{\lfloor T/M\rfloor M}).
\end{cases}
\]

Now $u_{\e}$ is an admissible function for the original problem, so the conclusion follows by \eqref{above4}; in particular
\[
F_{\e}(u_{\e},\Omega) = F_{\e}(u_{\e},\Omega\setminus B(z_{\e}, \e^{\lambda_1})) \leq \frac{C_{\rm hom}+o_M(1) + (M\log2)^{d-1} \omega(\delta/\e^{\lambda_1})}{(\lambda_1|\log\e|+\log d-\log2)^{d-1}}\,;
\]
hence
\[
\limsup_{\e\to0} |\log\e|^{d-1}\mu_{\e,\delta} \leq \inf_{\lambda_1 \in (0,1)}  C_{\rm hom}/\lambda_1^{d-1} = C_{\rm hom}. 
\]

\subsection{Proof of the main result about convergence of minima}

As a consequence of the previous section, we prove the main result on the asymptotic behaviour of minima defined in \eqref{minpbm1} by
\[
m_{\e,\delta}:=\min\Bigl\{\int_{\Omega} f\left(\frac{x}{\delta},\nabla u(x)\right)\,dx: u\in W^{1,d}_0(\Omega), u=1 \hbox{ on } B(z,\e), z\in\Omega\Bigr\},
\]
where also the centre of the small inclusion (a ball) is an argument of the minimization.

\begin{theorem} \label{mainpropd}
Assume there exists a point $x_0\in\Omega$ such that the following hold:

\emph{(i)} $f(x,\xi) \geq f(x_0,\xi)$ for every $x\in\mathbb{R}^d$ and for every $\xi\in\mathbb{R}^d$;

\emph{(ii)} for every $\nu>0$, there exists $r_{\nu}>0$ such that for every $x\in B(x_0, r_{\nu})$ and for every $\xi\in\mathbb{R}^d$ we have $f(x,\xi)\leq f(x_0,\xi) + \nu|\xi|^d$.

\noindent Then 
\[
\lim_{\e\to0}|\log\e|^{d-1}m_{\e,\delta}=
\Phi(x_0)C_{\rm hom}\Bigl[\lambda\Phi(x_0)^{\frac{1}{d-1}}+(1-\lambda)C_{\rm hom}^{\frac{1}{d-1}}\Bigr]^{1-d}.
\]
\end{theorem}

\begin{proof}
Since we use the same argument presented in the proof of Proposition \ref{mainpropfixed}, we focus on highlighting the main differences, keeping the same notations.

\smallbreak

\emph{Bound from below}. In the case $\lambda=0$, we introduce $\lambda_2>0$, then we apply Lemma \ref{modlemma} to get the inequality
\[
\left(1+\frac{C}{N-1}\right)F_\e(u,\Omega)\geq 
\]
\begin{align*}
\geq & \min\{F_{\e}(v,B(z,\e^{\lambda_2})) : v\in W^{1,d}_0(B(z,\e^{\lambda_2})), 
v=1 \hbox{ on } B(z,\e)\}|1-c|^d \\ \nonumber
+ & \min\{F_{\e}(v,B(z, R_{\Omega})\setminus  \overline{B}(z,\e2^{S-N})):v\in W^{1,d}(B(z, R_{\Omega})\setminus \overline{B}(z,\e^{\lambda_2}2^{S-N})), \\ \label{min2}
& \qquad \qquad \qquad \qquad \qquad \qquad v=1 \hbox{ on } B(z,\e2^{S-N}), v=0 \hbox{ on } \partial B(z, R_{\Omega})\}|c|^d.
\end{align*}
Note that the second summand is estimated exactly as \eqref{min2}; while, for the first summand, we cannot exploit the property of periodicity (P$1$) of the energy since minimization involves also the centre of the inclusion. To deal with this term, we consider a minimizer $u$ and we simply apply (i) to get
\begin{equation}\label{mainpropbelow}
\begin{split}
& \min\{F_{\e}(v, B(z, \e^{\lambda_2})) : v\in W^{1,d}_0(B(z,\e^{\lambda_2})), 
v=1 \hbox{ on } B(z,\e)\}|1-c|^d \\
& \qquad \qquad \qquad \qquad \geq \int_{B(z,\e^{\lambda_2})} f(x_0, \nabla u(x))\,dx \, |1-c|^d= \frac{\Phi(x_0)+o_\e(1)}{(1-\lambda_2)^{d-1}|\log \e|^{d-1}}|1-c|^d.
\end{split}
\end{equation}
This is the same estimate we obtained in \eqref{below2lambda0}, with the point $x_0$ in place of the fixed centre $z$. Analogously to Proposition \ref{mainpropfixed}, we conclude that $|\log\e|^{d-1}m_{\e,\delta}\to\Phi(x_0)$.

\smallbreak

If $\lambda\in(0,1)$, we further introduce $\lambda_1\in(0,\lambda)$ and we achieve the inequality
\[
\biggl(1+\frac{C}{N-1}\biggr)F_{\e}(u,\Omega) \geq 
\]
\begin{align}
& \geq \min\{F_{\e}(v,B(z_{\e},\e^{\lambda_2})) : v\in W^{1,d}_0(B(z_{\e},\e^{\lambda_2})), v=1 \hbox{ on } B(z_{\e},\e)\}|1-c_{-1}|^d \label{in01-1bis} \\
& + \min\{F_{\e}(v,B(z_{\e},\e^{\lambda_1})  \setminus B(z_{\e},\e2^{S'-N})) : v\in W^{1,d}(B(z_{\e},\e^{\lambda_1})\setminus \overline{B}(z_\e,\e2^{S'-N})), \notag \\ 
& \qquad \qquad \qquad \qquad v=1 \hbox{ on } \partial B(z_{\e},\e2^{S'-N}), v=0 \hbox{ on } \partial B(z_{\e},\e^{\lambda_1})\}|c_{-1}-c_0|^d \label{in01-2bis} \\
& + \min\{F_{\e}(v,B(z_{\e},\e^{\lambda_1}2^M) \setminus B(z_{\e},\e2^{S''-N})) : v\in W^{1,d}(B(z_{\e},\e^{\lambda_1}2^M) \setminus \overline{B}(z_{\e},\e2^{S''-N})), \notag \\ 
& \qquad \qquad \qquad \qquad v=1 \hbox{ on } \partial B(z_{\e},\e2^{S''-N}), v=0 \hbox{ on } \partial B(z_{\e},\e^{\lambda_1}2^M)\}|c_{0}-c_1|^d \label{in01-3bis} \\
& + \sum_{k=2}^{\lfloor T/M\rfloor+1} \min\{F_{\e}(v, A_{M,k}^N) : v\in W^{1,d}(A_{M,k}^N), \,  v=1 \hbox{ on } \partial B(z_{\e},\e^{\lambda_1}2^{(k-1)M-N}), \notag \\
& \qquad \qquad \qquad \qquad v=0 \hbox{ on } \partial B(z_{\e},\e^{\lambda_1}2^{kM})\}|c_{k-1}-c_{k}|^d\,, \label{in01-4bis}
\end{align}
where we put $c_{\lfloor \frac{T}{M} \rfloor+1}:=0$.

The estimates for the terms \eqref{in01-2bis}, \eqref{in01-3bis}, \eqref{in01-4bis}
are achieved precisely as in \eqref{in01-2}, \eqref{in01-3}, \eqref{in01-4} respectively, while \eqref{in01-1bis} is estimated exploiting (i) as in \eqref{mainpropbelow}. Once more, the outcome is the same of Proposition \ref{mainpropfixed}, with $x_0$ in place of $z$.

\smallbreak

The case $\lambda=1$ is analogous and can be proved starting by the estimate in \eqref{mainpropfixedhom}; this might be expected since, at this scale, the only effect in the minimization is due to the homogenization (and not to the point in which we concentrate our inclusion).

\smallbreak

\emph{Bound frome above}. Take $z_{\e}=\delta x_0$ modulo the $\delta$-cube in such a way that this family of points is contained in a ball $B\subset \subset \Omega$. Condition (ii) allows to apply the bound from above given by Proposition \ref{mainpropfixed}, then we conclude observing that $m_{\e,\delta} \leq \mu_{\e,\delta}$.
\end{proof}

\medbreak

We remark that assumption (i) may be weakened. Note indeed that the key estimate we need to carry out our proof, and more specifically the bound from below, is
\begin{multline*}
\min\{F_{\e}(v, B(z, \e^{\lambda_2})) : v\in W^{1,d}_0(B(z,\e^{\lambda_2})), 
v=1 \hbox{ on } B(z,\e)\} \\
\geq \int_{B(z,\e^{\lambda_2})} f(x_0, \nabla u(x))\,dx,
\end{multline*}
where $u$ is a minimizer for fixed $\lambda_2\in(\lambda,1)$.

\medbreak

A plausible sufficient condition might seem to be that $\Phi$ attains its minimum at the point $x_0$. Yet, note that this requirement is inadequate if $\Phi$ is not continuous at a minimum point. For instance, consider the function defined on $(0,1)^d$ as
\[
f(x,\xi) :=
\begin{cases}
\frac{1}{2}|\xi|^d & \hbox{ if } x=x_0:=\left(\frac{1}{2},...,\frac{1}{2}\right) \\[3pt]
|\xi|^d & \hbox{ otherwise} \\[3pt]
\end{cases}
\]
and then extended by periodicity; we see that \eqref{minpbm1} reduces to the homogeneous problem and then, that $|\log \e|^{d-1}m_{\e,\delta}\to \sigma_{d-1}$ as $\e\to0$. But this is a contradiction, indeed, we have $f_{\rm hom}(\xi)=|\xi|^d$ so that $C_{\rm hom}=\sigma_{d-1}$, while $\Phi(x_0)$ is equal to $\sigma_{d-1}/2$; plugging these in \eqref{result1} and assuming $\lambda=0$, we get $|\log \e|^{d-1}m_{\e,\delta}= \Phi(x_0)=\sigma_{d-1}/2$. 

\section{Application to perforated domains}

In this final section we maintain the setting and notation introduced in the previous ones. We will make use of Proposition \ref{mainpropfixed} to compute the $\Gamma$-limit of a family of functionals defined with boundary conditions related to varying domains.

We fix $(\e_k)_{k\in\mathbb{N}}$ a positive sequence converging to $0$, we consider the sequence of critical periods $d_{k}:=|\log \e_k|^{\frac{1-d}{d}}$, and for every $i\in\mathbb{Z}^d$, we put $x_{k}^i:=id_k$.

\noindent Then, we introduce a further scale which rules the periodic structure of the energy, say $\delta=\delta(\e)$, and we define $\delta_k:=\delta(\e_k)$ for every $k\in\mathbb{N}$, obtaining a positive sequence vanishing as $k\to+\infty$. 
In accordance with the previous sections, we will always assume that it exists 
\begin{equation}\label{lambda}
\lambda:=\lim_{k\to+\infty} \frac{|\log \delta_k|}{|\log \e_k|} \land 1.
\end{equation}
\noindent Assuming that $\Omega$ is a bounded open subset of $\mathbb{R}^d$ such that $|\partial \Omega|=0$, we define a periodically perforated domain as 
\[
\Omega_{k}:=\Omega \setminus \bigcup_{i\in \mathbb{Z}^d} B(x_{k}^i,\e_k),
\]
and we consider functionals $F_{k}: L^d(\Omega)\rightarrow[0,+\infty]$ given by
\[
F_{k}(u):=
\begin{cases}
\displaystyle \int_{\Omega}f\left(\frac{x}{\delta_k}, \nabla u(x)\right)\,dx & \hbox{ if } u\in W^{1,d}(\Omega) \hbox{ and } u=0 \hbox{ on } \Omega\setminus\Omega_k \\[5pt]
+\infty & \hbox{ otherwise}.
\end{cases}
\]

To prove our result, we assume that the perforations are related to the periodic structure of the heterogeneous medium, in particular we suppose that 
\begin{equation}\label{periodicity}
\hbox{for every } k \hbox{ there exists a positive natural number } m_k \hbox{ such that } d_k=m_k\delta_k 
\end{equation}
and that
\begin{equation}\label{hyplimit}
\frac{\delta_k}{d_k}\to0 \hbox{ as } k\to+\infty.
\end{equation}
\noindent Condition \eqref{periodicity} leads to the identity
\begin{equation}\label{hypperiod}
f\left(\frac{x_k^i}{\delta_k}+y,\xi\right)=f\left(\frac{id_k}{\delta_k}+y,\xi\right)=f(y,\xi) \hbox{ for every } i\in\mathbb{Z}^d,\,y\in\mathbb{R}^d, \xi\in\mathbb{R}^d.
\end{equation}
If \eqref{periodicity} is not fulfilled, then $f\left(\frac{x_k^i}{\delta_k}+y,\xi\right)=f(y_k+y,\xi)$ for some $y_k\in [0,1]^d$, and the result depends on the properties of $(y_k)_k$ modulo $\delta$, see \cite{AnsBra1} for the occurrence of a similar phenomenon.

\smallbreak

In order to apply Proposition \ref{mainpropfixed}, we add suitable regularity assumptions on $f$ at the point $0$. Our statement reads as follows.

\begin{theorem}\label{gammalimit}
Assume that for every $\nu>0$, there exists $r_{\nu}>0$ such that for every $x\in B(0, r_{\nu})$ it holds 
\begin{equation}\label{hyp}
|f(0,\xi)-f(x,\xi)|\leq \nu|\xi|^d \hbox{ for every } \xi\in\mathbb{R}^d.
\end{equation}

\noindent Then
\[
\Gamma\emph{-}\lim_k F_k(u)=F(u):=\int_{\Omega}f_{\rm hom}(\nabla u(x))\,dx + C(\lambda)\int_{\Omega}|u(x)|^d\,dx,
\]
for every $u\in W^{1,d}(\Omega)$, where the $\Gamma$-limit is meant with respect to the strong convergence in $L^d(\Omega)$ and $C(\lambda)$ is given by
\[
C(\lambda):=
\Phi(0)C_{\rm hom}\Bigl[\lambda\Phi(0)^{\frac{1}{d-1}}+(1-\lambda)C_{\rm hom}^{\frac{1}{d-1}}\Bigr]^{1-d},
\]
with $\Phi, C_{\rm hom}, f_{\rm hom}$ and $\lambda$ defined as in \eqref{Phi}, \eqref{Chom}, \eqref{fhom} and \eqref{lambda} respectively.

\end{theorem}

We basically prove that, in the $\Gamma$-limit, internal boundary conditions imposed on the perforations vanish, being replaced by the additional term $C(\lambda)\int_{\Omega}|u|^ddx$.

\subsection{The main construction and some auxiliary results}

In our proof we will make wide use of Lemma \ref{modlemma}, but its application is more delicate in this instance. To fit our arguments, it needs some refinement: we perform the modifications among annuli which are not only homothetic, but also such that their corresponding inner and outer radii are proportional to the period $d_k$.

\smallbreak

We introduce $Z_{k}:=\{i\in\mathbb{Z}^d : \text{dist}(x_k^i,\partial\Omega)>d_k\}$, namely the set of the centres of those perforations which are uniformly far from the boundary.

Let $M\in\mathbb{N}$, $\alpha>0$ be such that $\alpha2^{M+1}<1/2$. Given a sequence $(u_k)_k$ in $W^{1,d}(\Omega)$, fix $k$, and around each point $x_k^i$ with $i\in Z_k$ apply Lemma \ref{modlemma} to the function $u_k$ with
\begin{equation}\label{applemma}
f(x,\xi)=f\left(\frac{x}{\delta},\xi\right)\,,\eta=\alpha d_k\,, R=\alpha2^{M+1}d_k\,, N=M \hbox{ and } r=\alpha d_k .
\end{equation}

We obtain a function $v_k$ having constant values $u_k^i$ on the boundary of each ball centered at $x_k^i$ with radius $\alpha 2^{j_i}d_k$ for some $j_i\in\{1,...,M\}$ and $i\in Z_k$. Also recall that this function comes with the estimate
\[
\int_{\Omega}f\left(\frac{x}{\delta_k}, \nabla v_k(x)\right)\,dx \leq \left(1+\frac{C}{M-1}\right) \int_{\Omega}f\left(\frac{x}{\delta_k}, \nabla u_k(x)\right)\,dx\,.
\]

\medbreak

We take advantage of the following result which is a simplified version of the discretization argument proved by Sigalotti (see \cite[Proposition 3.3]{Sigalotti}).

\begin{proposition}\label{cubes}
Let $(u_k)_k$ be a sequence in $W^{1,d}(\Omega)\cap L^{\infty}(\Omega)$ strongly converging to $u$ in $L^d(\Omega)$ and such that $(\nabla u_k)_k \subseteq L^d(\Omega)$ is bounded. For every $i\in Z_k$, let $u_k^i$ be the mean values described above and put 
\[
Q_k^i:=x_k^i+\left(-\frac{d_k}{2}, \frac{d_k}{2}\right)^d.
\]

Then
\[
\lim_{k\to\infty}\int_{\Omega}\Biggl|\sum_{i\in Z_k} |u_k^i|^d\chi_{Q_k^i}(x)-|u(x)|^d\Biggr|\,dx = 0.
\]
\end{proposition}

A useful tool to proceed will be also the following convergence result which is an application of the Riemann-Lebesgue lemma.

\begin{lemma}\label{riemannlebesgue} The sequence 
\[
\chi_k(x):=\chi_{\Omega \setminus \bigcup_{i\in Z_k}B(x_k^i,d_k/2)}(x)\,, \qquad k\in\mathbb{N}
\]
weakly* converges to a positive constant in $L^{\infty}(\Omega)$.
\end{lemma}

\subsection{Liminf inequality}

We prove that for every $u\in W^{1,d}(\Omega)$ and for every sequence $(u_k)_k$ in $L^d(\Omega)$ such that $u_k \to u$ in $L^d(\Omega)$, it holds $\liminf_kF_k(u_k) \geq F(u)$. 

\medbreak

The first step of the proof consists in applying the modification lemma as in \eqref{applemma}. To simplify the notation in this section, we limit ourselves to denote the radii on which the modified function $v_k$ attains the constant values $u_k^i$ by $\rho_k^i$ in place of $\alpha 2^{j_i}d_k$.

\smallbreak

Without loss of generality we may assume $(u_k)_k\subseteq W^{1,d}(\Omega)$ and $\sup_kF_k(u_k)<+\infty$. Note that the last condition, combined with the equi-coerciveness of the functionals $(F_k)_k$, implies that $\sup_k\|\nabla u_k\|_{L^d(\Omega)}<\infty$, hence $u_k \rightharpoonup u$ in $W^{1,d}(\Omega)$.

\smallbreak

In a first instance, also assume that $(u_k)_k$ is bounded in $L^{\infty}(\Omega)$. We aim to estimate
\begin{equation} \label{inliminf1}
F_k(v_k)=\int_{\Omega \setminus \bigcup_{i\in Z_k} B(x_k^i, \rho_k^i)}f\left(\frac{x}{\delta_k},\nabla v_k(x)\right)\,dx + \sum_{i \in Z_k}\int_{B(x_k^i, \rho_k^i)} f\left(\frac{x}{\delta_k}, \nabla v_k(x)\right)\,dx\,.
\end{equation}

We perform another modification putting
\[
w_k:=
\begin{cases}
v_k & \hbox{ on } \Omega \setminus \bigcup_{i\in Z_k} B(x_k^i, \rho_k^i), \\
u_k^i & \hbox{ on } B(x_k^i, \rho_k^i),\, i\in Z_k.
\end{cases}
\]

It trivially holds
\[
\int_{\Omega \setminus \bigcup_{i\in Z_k} B(x_k^i, \rho_k^i)}f\left(\frac{x}{\delta_k},\nabla v_k(x)\right)\,dx = \int_{\Omega} f\left(\frac{x}{\delta_k},\nabla w_k(x)\right)\,dx\,.
\]

Note that, according to the proof of Lemma \ref{modlemma}, $\|v_k\|_{L^{\infty}(\Omega)} \leq  \|u_k\|_{L^{\infty}(\Omega)}$, hence $\|w_k\|_{L^{\infty}(\Omega)} \leq  \|u_k\|_{L^{\infty}(\Omega)}$ so that $(w_k)_k$ is bounded in $L^{\infty}(\Omega)$ and then also bounded in $L^d(\Omega)$. Moreover, as $\left(1+\frac{C}{M-1}\right)F_k(u_k) \geq F_k(v_k) \geq F_k(w_k)$, we deduce that $(w_k)_k$ is bounded in $W^{1,d}(\Omega)$; thus, we may extract a subsequence $(w_{k_j})_j$ weakly converging to a certain $w$ in $W^{1,d}(\Omega)$. 

\noindent As $w_k-u_k \in W^{1,d}_0(\Omega)$ for every $k$ and since $u_k\rightharpoonup u$ in $W^{1,d}(\Omega)$ and $u_k \to u$ in $L^d(\Omega)$, it holds by Rellich's Theorem that $(w_{k_j})_j$ actually converges strongly to $w$ in $L^d(\Omega)$. We claim that such $w$ does not depend on the subsequence and that it coincides with $u$.
\noindent To prove this, note that for every $k$
\[
w_k\chi_{\Omega \setminus \bigcup_{i\in Z_k}B(x_k^i,d_k/2)}=u_k\chi_{\Omega \setminus \bigcup_{i\in Z_k}B(x_k^i,d_k/2)}
\]
and also that by Lemma \ref{riemannlebesgue} and the previous observations, the following hold
\[
\begin{cases}
\chi_{\Omega \setminus \bigcup_{i\in Z_k}B(x_k^i,d_k/2)} \overset{*}{\rightharpoonup} c & \hbox{ in }
L^{\infty}(\Omega), \\
u_k\to u & \hbox{ in } L^d(\Omega), \\
w_{k_j} \to w & \hbox{ in } L^d(\Omega).
\end{cases}
\]
These facts imply
\[
\begin{cases}
\chi_{\Omega \setminus \bigcup_{i\in Z_k}B(x_k^i,d_k/2)}u_k \rightharpoonup cu & \hbox{ in } L^d(\Omega), \\
\chi_{\Omega \setminus \bigcup_{i\in Z_{k_j}}B(x_{k_j}^i,d_{k_j}/2)}w_{k_j} \rightharpoonup cw & \hbox{ in } L^d(\Omega),
\end{cases}
\]
hence, it follows that $u=w$ in $L^d(\Omega)$ for every subsequence, proving that $w_k \to u$ in $L^d(\Omega)$.
By the Homogenization Theorem and the liminf inequality, we deduce
\begin{equation} \label{liminf1}
\liminf_k \int_{\Omega} f\left(\frac{x}{\delta_k},\nabla w_k(x)\right)\,dx \geq \int_{\Omega}f_{\rm hom}(\nabla u(x))\,dx\,.    
\end{equation}

To estimate the second contribution in \eqref{inliminf1}, fix $i\in Z_k$ and let $\varphi_k^i$ be a function solving
\[
\min \Bigl\{\int_{B(x_k^i,\rho_k^i)} f\left(\frac{x}{\delta_k},\nabla u(x)\right)\,dx : u\in u_k^i+W^{1,d}_0(B(x_k^i, \rho_k^i)), u=0 \hbox{ on } B(x_k^i, \e_k)\Bigr\}.
\]
Up to extending the function $\varphi_k^i$ to the constant $u_k^i$ on $B(x_k^i, d_k/2)\setminus B(x_k^i, \rho_k^i)$, we have
\[
\int_{B(x_k^i, \rho_k^i)} f\left(\frac{x}{\delta_k},\nabla v_k(x)\right)\,dx 
\geq \int_{B(x_k^i, \rho_k^i)} f\left(\frac{x}{\delta_k},\nabla \varphi_k^i(x)\right)dx 
\]
\begin{equation*} 
\begin{split}
    & \geq \min \Bigl\{\int_{(B\left(x_k^i,\frac{d_k}{2}\right)} f\left(\frac{x}{\delta_k},\nabla u\right)dx : u\in u_k^i+W^{1,d}_0(B(x_k^i, d_k/2)), u=0 \hbox{ on } B(x_k^i, \e_k)\Bigr\} \\
& = \min\Bigl\{\int_{B\left(0,\frac{1}{2}\right)} f\left(\frac{d_k x}{\delta_k},\nabla u\right)dx : u\in 1+W^{1,d}_0(B(0,1/2), u=0 \hbox{ on } B(0,\e_k/d_k)\Bigr\}|u_k^i|^d,
\end{split}
\end{equation*}
where the last equality follows by the change of variables $x\mapsto x_k^i+d_kx$, the identity \eqref{hypperiod} and (P$2$).

\noindent Now put
\[
\delta'_k:=\frac{\delta_k}{d_k}\,, \qquad \e'_k:=\frac{\e_k}{d_k}\,, \qquad  \lambda':=\lim_k \frac{|\log \delta'_k|}{|\log \e'_k|},
\] 
and rewrite the previous inequality as
\begin{multline}\label{liminfineq1}
\int_{B(x_k^i, \rho_k^i)} f\left(\frac{x}{\delta_k},\nabla v_k(x)\right)\,dx \\
\geq \min\Bigl\{\int_{B\left(0,\frac{1}{2}\right)} f\left(\frac{x}{\delta'_k},\nabla u\right)dx : u\in 1+W^{1,d}_0(B(0,1/2), u=0 \hbox{ on } B(0,\e'_k)\Bigr\}|u_k^i|^d.
\end{multline}

Note that $\e'_k= \e_k|\log\e_k|^{1-1/d}\to0$, while $\delta'_k=\delta_k|\log\e_k|^{1-1/d}\to0$ as $k\to\infty$ by assumption \eqref{hyplimit}; also observe that
\[
\lambda'=\lim_k \frac{|\log\delta_k+\log |\log\e_k|^{1-1/d}|}{|\log \e_k+\log|\log \e_k|^{1-1/d}|}=\lim_k\frac{|\log\delta_k|}{|\log \e_k|}=\lambda.
\]
In light of the assumption \eqref{hyp}, we are in position to apply Proposition \ref{mainpropfixed} (up to the transformation $u\mapsto 1-u$) to \eqref{liminfineq1} with $\Omega=B(0,1/2)$ and $z_{\e}=0$ for every $\e$. We get
\[
\min\Bigl\{\int_{B(0,1/2)} f\left(\frac{x}{\delta'_k},\nabla u(x)\right)\,dx : u\in 1+W^{1,d}_0(B(0, 1/2)), u=0 \hbox{ on } B(0,\e'_k)\Bigr\} =
\]
\[
= \frac{C(\lambda)+o_k(1)}{|\log\e'_k|^{d-1}} = \frac{C(\lambda)+o_k(1)}{|\log\e_k|^{d-1}}\,,
\]
and by Proposition \ref{cubes}, it follows
\begin{equation}\label{liminf2}
\begin{split}
\liminf_k \sum_{i \in Z_k}\int_{B(x_k^i, d_k/2)} f\left(\frac{x}{\delta_k}, \nabla v_k(x)\right)\,dx & \geq \liminf_k \frac{C(\lambda)}{|\log\e_k|^{d-1}}\sum_{i \in Z_k}|u_k^i|^d+o_k(1) \\
& = C(\lambda) \int_{\Omega}|u(x)|^d\,dx\,.
\end{split}
\end{equation}

Finally, by \eqref{liminf1} and \eqref{liminf2}, we deduce
\[
\left(1+\frac{C}{M-1}\right)\liminf_k F_k(u_k) \geq \liminf_k F_k(v_k) \geq \int_{\Omega}f_{\rm hom}(\nabla u(x))\,dx + C(\lambda)\int_{\Omega}|u(x)|^d\,dx.
\]

Recall that $\alpha$ and $M$ have been chosen so that $\alpha 2^{M+1} < 1/2$ and, since the reasoning leading to the above estimate holds true for every $\alpha>0$, we may let $M\to+\infty$ getting the liminf inequality.

\medbreak

We conclude removing the boundedness assumption on $(u_k)_k\subseteq L^{\infty}(\Omega)$ by a truncation argument: assume $u_k\to u$ in $L^d(\Omega)$ and put $\overline{u}_k^M:=((-M) \lor u_k)\land M$ for fixed $M\in\mathbb{N}$; by dominated convergence, $\overline{u}_k^M\to u$ in $L^d(\Omega)$ as $k,M\to+\infty$, moreover, since $f(\cdot,0)=0$, it holds 
\[
\int_{\Omega}f\left(\frac{x}{\delta_k}, \nabla u_k\right)\,dx \geq \int_{\Omega}f\left(\frac{x}{\delta_k}, \nabla \overline{u}_k^M\right)\,dx
\]
for every $k,M \in\mathbb{N}$, thus we immediately conclude by the previous instance.

\medbreak

Denoting by $F':=\Gamma$-$\liminf_k F_k$, what we have proved so far is that $F(u)\leq F'(u)$ for every $u\in W^{1,d}(\Omega)$.

\subsection{Limsup inequality}

The goal of this section is to define a recovery sequence converging in $L^d(\Omega)$ to a fixed function $u\in W^{1,d}(\Omega)$. First we assume that $u\in L^{\infty}(\Omega)$.

Start by a recovery sequence $u_k\to u$ in $L^d(\Omega)$ related to the functionals
\[
F_k^0(u):=
\begin{cases}
\displaystyle \int_{\Omega}f\left(\frac{x}{\delta_k}, \nabla u(x)\right)\,dx & \hbox{ if } u\in W^{1,d}(\Omega), \\[5pt]
+\infty & \hbox{ if } u\in L^d(\Omega)\setminus W^{1,d}(\Omega)
\end{cases}
\]
which are known to $\Gamma$-converge to 
\[
F^0(u):=\int_{\Omega}f_{\rm hom}(\nabla u(x))\,dx
\]
for every $u\in W^{1,d}(\Omega)$ as stated in the Homogenization Theorem. By the equi-coerciveness of the functionals $(F^0_k)_k$, we deduce $u_k \rightharpoonup u$ in $W^{1,d}(\Omega)$. 

\noindent It is a known fact that, up to extract a subsequence, we can also assume that $(|\nabla u_k|^d)_k$ is an equi-integrable family (see \cite{FMP} and \cite[Remark C.6]{BDF}).

We claim that we can also make our recovery sequence bounded in $L^{\infty}(\Omega)$.
Let $T:=\|u\|_{L^{\infty}(\Omega)}$ and define $u_k':=(-(T+1) \lor u_k )\land (T+1)$. We get a bounded sequence in $L^{\infty}(\Omega)$ which converges to $u$ in $L^d(\Omega)$ with further property that $(|\nabla u'_k|^d)_k$ is still equi-integrable being obtained by truncation.

\noindent Note that 
\begin{multline*}
\left|\int_{\Omega}f_{\rm hom}(\nabla u_k(x))\,dx - \int_{\Omega}f_{\rm hom}(\nabla u'_k(x))\,dx\right|  \leq \int_{\{|u_k|>T+1\}}|f_{\rm hom}(\nabla u_k(x))|\,dx \\
\leq \beta \int_{\{|u_k|>T+1\}}|\nabla u_k(x))|^d\,dx \leq \beta \int_{\{|u_k-u|>1\}}|\nabla u_k(x))|^d\,dx\,;
\end{multline*}
but since $u_k\to u$ in measure and $(|\nabla u_k|^d)_k$ is equi-intergable, the last term tends to $0$ and the claim is proved.

\smallbreak

For every $k$, define modifications $v_k$ by transformations around every point $x_k^i$ with $i\in Z_k$ as we did in \eqref{applemma}. We recall the construction for clarity: fix $M\in\mathbb{N}$ and let $\alpha>0$ be such that $\alpha2^{M+1}<1/2$, then apply Lemma \ref{modlemma} with
\[
f(x,\xi)=f\left(\frac{x}{\delta},\xi\right)\,,\eta=\alpha d_k\,, R=\alpha2^{M+1}d_k\,, N=M \hbox{ and } r=\alpha d_k .
\]
We have that
\[
\int_{\Omega}f\left(\frac{x}{\delta_k}, \nabla v_k(x)\right)\,dx \leq \left(1+\frac{C}{M-1}\right) \int_{\Omega}f\left(\frac{x}{\delta_k}, \nabla u_k(x)\right)\,dx\,,
\]
and the function $v_k$ attains the constant value $u_k^i$ on $\partial B(x_k^i,\rho_k^i)$, where $\rho_k^i$ is of the form $\alpha2^{j_i}d_k$ for some $j_i\in\{1,...,M\}$.

\noindent Since $\e_k/d_k\to0$ as $k\to+\infty$, we can also assume $\e_k<\alpha d_k$ for every $k$; hence, we define
\[
w_k:=
\begin{cases}
v_k & \hbox{ on } \Omega\setminus \bigcup_{i\in Z_k} B(x_k^i,\rho_k^i)\\
u_k^i & \hbox{ on } B(x_k^i, \rho_k^i)\setminus B(x_k^i, \alpha d_k), i\in Z_k \\
\varphi_k^i & \hbox{ on }  B(x_k^i, \alpha d_k),\, i\in Z_k,
\end{cases}
\]
where $\varphi_k^i$ solves the minimum problem 
\begin{multline*}
\min\Bigl\{\int_{B(x_k^i, \alpha d_k)}f\left(\frac{x}{\delta_k}, \nabla u(x)\right)\,dx : u \in u_k^i+W^{1,d}_0(B(x_k^i, \alpha d_k)), u = 0 \hbox{ on } B(x_k^i, \e_k)\Bigr\}\\
= \min\Bigl\{\int_{B(0, \alpha)}f\left(\frac{x}{\delta'_k}, \nabla u(x)\right)\,dx : u \in 1+W^{1,d}_0(B(0, \alpha)), u = 0 \hbox{ on } B(0, \e'_k)\Bigr\}|u_k^i|^d \\
= \frac{C(\lambda)+o_k(1)}{|\log\e_k|^{d-1}}|u_k^i|^d\,.
\end{multline*}

Let $A_k:= \bigcup_{i\in Z_k}B(x_k^i, \rho_k^i)$. We will treat with different arguments the contributions due to $\Omega\setminus A_k$ and $A_k$. 

\smallbreak

We estimate the contribution on $A_k$ using Proposition \ref{cubes},
\begin{multline}\label{limsup1}
\limsup_k \int_{A_k}f\left(\frac{x}{\delta_k}, \nabla w_k(x)\right)\,dx = \limsup_k \sum_{i\in Z_k}  \int_{B(x_k^i,\alpha d_k)}f\left(\frac{x}{\delta_k},\nabla \varphi_k^i(x)\right)\,dx \\
= \limsup_k \sum_{i\in Z_k} |u_k^i|^d \frac{C(\lambda)+o_k(1)}{|\log\e_k|^{d-1}} = C(\lambda)\int_{\Omega}|u(x)|^d\,dx\,.
\end{multline}

To estimate the contribution on $\Omega\setminus A_k$, we put 
\[
Z_k':=\{i\in\mathbb{Z}^d : B(x_k^i, \e_k) \cap \Omega \neq \emptyset, i\notin Z_k\}\qquad , \qquad r_k:=\alpha2^{M+1}d_k 
\]
and define 
\[
A'_k:=\bigcup_{i\in Z'_k}B(x_k^i,r_k)
\]
in order to study separately the behaviours on $A'_k$ and $\Omega \setminus (A_k \cup A'_k)$.

\smallbreak

Take into account the contribution of $A'_k$. We set
\[
Q_k^i:=x_k^i+\left(-\frac{d_k}{2}, \frac{d_k}{2}\right),
\]
and we see preliminarily that
\begin{equation}\label{boundary}
|\Omega \cap A'_k| \leq \sum_{i\in Z'_k} (r_k)^d \sim \# Z'_k(d_k)^d =\biggl| \bigcup_{i\in Z'_k} Q_k^i\biggr|\to|\partial \Omega|=0
\end{equation}
by assumption.

For every $i\in Z'_k$, let $\psi_k^i$ be the solution to the homogeneous capacitary problem
\[
\min\Bigl\{\int_{B(x_k^i, r_k)} |\nabla u(x)|^d\,dx : u \in 1+W^{1,d}_0(B(x_k^i,r_k)), u = 0 \hbox{ on } B(x_k^i, \e_k)\Bigr\}
\]
which is (known to be) equal to $\sigma_{d-1}|\log r_k-\log\e_k|^{1-d}$. 

\noindent Up to extending $\psi_k^i$ with value $1$ on $\mathbb{R}^d\setminus B(x_k^i,r_k)$, we set as recovery sequence 
\[
w'_k:=w_k \prod_{i\in Z'_k}\psi_k^i \qquad \hbox{ on } \Omega.
\]

Such $w'_k$ is a modification of $w_k$ performed on $A'_k$, which is disjoint from $A_k$ in virtue of the choice of the radii $r_k$ and $\rho_k^i$, for every $k$ and $i\in Z_k$. This means that the estimate in \eqref{limsup1} is still valid replacing $w_k$ with $w'_k$.

\smallbreak

We prove that 
\begin{equation}\label{limsup2}
\limsup_k \int_{\Omega \cap A'_k}f\left(\frac{x}{\delta_k}, \nabla w'_k(x)\right)\,dx = 0.
\end{equation}
For every $i\in Z'_k$, we have 
\begin{multline*}
\int_{\Omega \cap B(x_k^i,r_k)} f\left(\frac{x}{\delta_k}, \nabla w'_k(x)\right) 
\leq \beta \int_{\Omega \cap B(x_k^i,r_k)}|\nabla w'_k(x)|^d\,dx \\
\leq 2^{d-1}\beta\left[(1+ \|u\|_{L^{\infty}(\Omega)})^d\int_{B(x_k^i,r_k)}|\nabla \psi_k^i(x)|^d\,dx + \int_{\Omega \cap B(x_k^i,r_k)} |\nabla w_k(x)|^d\,dx\right] \\
\leq C \left[|\log r_k-\log\e_k|^{1-d} + \int_{\Omega \cap B(x_k^i,r_k)} |\nabla w_k(x)|^d\,dx\right]
\end{multline*}
for a positive constant $C$ which depends only on $\|u\|_{L^{\infty}(\Omega)},\beta$ and the dimension $d$.

\noindent Note that, since $i\in Z'_k$, by definition of $w_k$ we have 
\[
\int_{\Omega \cap B(x_k^i,r_k)} |\nabla w_k(x)|^d\,dx = \int_{\Omega \cap B(x_k^i,r_k)} |\nabla v_k(x)|^d\,dx\,,
\]
and by the property (ii) of Lemma \ref{modlemma}, i.e., modifications on the starting function occur very close to the prescribed radius, it also holds
\[
\int_{\Omega \cap B(x_k^i,r_k)} |\nabla v_k(x)|^d\,dx = \int_{\Omega \cap B(x_k^i,r_k)} |\nabla u_k(x)|^d\,dx\,.
\]
Exploiting the equi-integrability of $(|\nabla u_k|^d)_k$, by \eqref{boundary} we infer that 
\[
\limsup_k \sum_{i\in Z'_k} \int_{\Omega \cap B(x_k^i,r_k)} |\nabla w_k(x)|^d\,dx = 0.
\]
At this point
\[
\limsup_k \int_{\Omega \cap A'_k}f\left(\frac{x}{\delta_k}, \nabla w'_k(x)\right)\,dx \leq C \limsup_k \sum_{i\in Z'_k} |\log r_k-\log\e_k|^{1-d}, 
\]
but since $\e_k \ll d_k $, we conclude that
\[
\limsup_k \sum_{i\in Z'_k} |\log\e_k|^{1-d} = \limsup_k \# Z'_k (d_k)^d = 0
\]
again by \eqref{boundary}.

\smallbreak

Finally, we deal with the contribution on $\Omega \setminus (A_k \cup A'_k)$. It holds
\begin{multline}\label{limsup3}
\limsup_k \int_{\Omega\setminus (A_k \cup A'_k)}f\left(\frac{x}{\delta_k}, \nabla w'_k(x)\right)\,dx = \limsup_k \int_{\Omega\setminus (A_k \cup A'_k)}f\left(\frac{x}{\delta_k}, \nabla v_k(x)\right)\,dx \\
\leq \limsup_k \int_{\Omega}f\left(\frac{x}{\delta_k}, \nabla v_k(x)\right)\,dx 
\leq \left(1+\frac{C}{M-1}\right)\limsup_k \int_{\Omega}f\left(\frac{x}{\delta_k}, \nabla u_k(x)\right)\,dx \\ 
\leq \left(1+\frac{C}{M-1}\right) \int_{\Omega}f_{\rm hom}(\nabla u(x))\,dx\,,
\end{multline}
where the last inequality is due to the fact that $(u_k)_k$ was originally picked as a recovery sequence to $u$ for the functionals $(F^0_k)_k$.

\smallbreak

Gathering \eqref{limsup1}, \eqref{limsup2} and \eqref{limsup3}, we get
\[
\limsup_k \int_{\Omega}f\left(\frac{x}{\delta_k}, \nabla w'_k\right)dx  \leq \left(1+\frac{C}{M-1}\right) \int_{\Omega}f_{\rm hom}(\nabla u)dx  + C(\lambda)\int_{\Omega}|u|^ddx\,.
\]
Since we can repeat the argument for every $\alpha>0$, we are free to set $M$ arbitrarily large, thus, the approximate limsup inequality is proved.

\medbreak

We still have to check that $w'_k\to u$ in $L^d(\Omega)$, i.e., it is actually an (approximate) recovery sequence.

\noindent Note that $\lim_k|\{w'_k \neq w_k\}|=0$ and $\sup_k\|w'_k-w_k\|_{L^{\infty}(\Omega)}\leq \|u_k\|_{L^{\infty}(\Omega)} \leq 1+\|u\|_{L^{\infty}(\Omega)}$ imply that $w'_k-w_k\to0$ in $L^d(\Omega)$, hence, it suffices to prove that $w_k\to u$ in $L^d(\Omega)$.

\noindent Since $\lim_k|\{w_k \neq v_k\}|=0$ and $\sup_k\|w_k-v_k\|_{L^{\infty}(\Omega)}\leq \|u_k\|_{L^{\infty}(\Omega)} \leq 1+\|u\|_{L^{\infty}(\Omega)}$, it holds that $w_k-v_k\to0$ in $L^d(\Omega)$, moreover $v_k\to u$ in $L^d(\Omega)$ by the same argument we used in the proof of the liminf inequality based on Lemma \ref{riemannlebesgue}; hence, $w_k\to u$ in $L^d(\Omega)$.

\bigbreak

To conclude, we remove the assumption $u\in L^{\infty}(\Omega)$. Recall that the $\Gamma$-limsup of $(F_k)_k$ is defined as
\[
F''(u):=\inf\{\hbox{lim\,sup}_k\, F_k(u_k) : u_k \to u \in L^d(\Omega)\}
\]
for every $u\in W^{1,d}(\Omega)$. $F''$ is sequentially lower semicontinuous with respect to the strong convergence in $L^d(\Omega)$ and by what we have already shown, it coincides with $F$ on $W^{1,d}(\Omega)\cap L^{\infty}(\Omega)$. 

\noindent Hence, given a sequence $(u_k)_k\subseteq W^{1,d}(\Omega)\cap L^{\infty}(\Omega)$ converging to $u$ in $ W^{1,d}(\Omega)$, it holds
\[
F''(u) \leq \liminf_k F''(u_k) = \liminf_k F(u_k)=F(u)
\]
by the continuity of $F$ with respect to the strong convergence in $W^{1,d}(\Omega)$, and this concludes the proof of the $\Gamma$-convergence.

\end{document}